\documentclass[10pt,  a4paper, reqno]{amsart}   
\usepackage{graphicx, caption, subcaption}
\usepackage{amsmath,mathtools,enumerate}
 
\usepackage{geometry}        

\usepackage[parfill]{parskip}  
\usepackage{graphicx}
\usepackage{textcomp} 
\usepackage[normalem]{ulem} 

\usepackage[applemac]{inputenc}
\usepackage{ae, aeguill}

\usepackage{pdfsync}

\usepackage{manfnt} 

\usepackage{pifont}  
\usepackage{qtree} 

\newcommand{\cleqn}{\setcounter{equation}{0}} 
\cleqn 

\newtheorem{theorem}{Theorem}[section]
\newtheorem{lemma}[theorem]{Lemma}
\newtheorem{prop}[theorem]{Proposition}

\everymath{\displaystyle}

\def \N{{\mathbb N}}
\def \Q{{\mathbb Q}}
\def \R{{\mathbb R}} 
\def \Z{{\mathbb Z}}

\def\A{{\mathcal A}}

\def\C{{\mathcal{ C}}}

\def\EE{{\mathbb{E}}}

\def\F{{\mathcal F}}

\def\K{{\mathcal K}}
\def\L{{\mathcal{L} }}

\def\M{{\mathcal M}}

\def\O{{\Omega}}
\def\cO{{\mathcal O}}

\def\PP{{\mathbb{ P}}}

\def\S{{\mathcal{S} }}

\def\tV{{\widetilde V}}
\def\uV{{\underline V}}

\def\tW{{\widetilde W}}

\def\a{\alpha}
\def\ta{\widetilde\alpha}
\def\o{\omega}  
\def\l{\lambda} 
\def\g{\gamma}

\def\d{\delta}

\def\1{{\mathbf 1}}

\def\Cov{\mathop{\mathbb{C}\rm ov}}

\def\eps{\varepsilon}

\def\Leb{\mathop{\rm Leb}}

\def\oB{\overline B}
\def\oH{\overline H}

\def\oL{\overline L}
\def\om{\overline m}

\def\phii{\varphi}

\def\supp{\mathop{\rm supp}}

\def\conl2{\stackrel $L^2$ \rightarrow}

\def\conl{\stackrel {\mathcal L} \rightarrow}

\def\bi{\bigbreak}

\def\qed{\hfill\rule{.2cm}{.2cm}}

\begin{document}

\title{\textbf{Quenched large deviations for Brownian motion in  a random potential}}
\author{Daniel Boivin and Thi Thu Hien L\^e  }
\date{\today}

\begin{abstract}
A quenched large deviation principle for  Brownian motion in
a non-negative,  stationary   potential   is proved. 
A sufficient moment condition on the potential is given but unlike
the results of Armstrong and Tran  (2014)  no regularity is assumed. 
The proof is  based on a method developed by Sznitman (1994) for
  Brownian motion among Poissonian potential. In particular, the LDP holds for
 potentials with polynomially decaying correlations such as the classical potentials
 studied by L. Pastur (1977) and R. Fukushima (2008)    
 and  the potentials recently introduced by H. Lacoin (2012).

\smallbreak\noindent
\end{abstract}

\maketitle
\markboth{D.Boivin and Hien L\^e}{Large deviations in a stationary potential}

\thanks{}


\section{Introduction}
\cleqn

Consider a standard Brownian motion on $\mathbb{R}^d$, $(Z_s; s\geq 0)$, moving in a non-negative
stationary  ergodic potential. That is, it is assumed that the potential is of the form
\begin{align*}
V(x, \o) := V_0(\tau_x \o ), \quad x\in\R^d,\o\in \O
\end{align*}
where $V_0$ is a real-valued non-negative random variable not identically zero on a probability space $(\O, \F, \PP)$
and  $({\tau}_{x}; x\in \R^d)$ is a    family of measurable maps on  $(\Omega,\F, \PP )$  which verifies
  \begin{align*}
  & \tau_{x}\circ \tau_{y}=\tau_{x+y}\  \hbox{ for all } x,y \in \R^d,\\
  & (x,\o)\mapsto \tau_x\o\quad\hbox{is measurable on the cartesian product  }\ \R^d\times \O, \\
  & \PP\  \hbox{is invariant under  }\   \tau_x \hbox{  for all  } x\in\R^d\hbox{  and}\\
  &  \hbox{is ergodic, that is, if for some  }  A\in\F,  \tau_x (A) = A  \hbox{  for all  }   x\in \mathbb{R}^d  \hbox{ then } \PP(A) = 0\hbox{ or } 1.
  \end{align*}
  
The \textit{quenched}  \textit{path measures}   are defined  by
\begin{equation} \label{quenched path}
 Q_{t,\omega}:=\frac{1}{S_{t,\omega}} \exp \big( -\int _{0}^{t} V(Z_s,\omega)ds \big) P_{0}, \hspace{0.5cm} t> 0, \omega \in \Omega
\end{equation}
where the normalizing constants $S_{t,\o}$ are the   \textit{quenched  survival functions up to time}  $t$
\begin{equation}\label{quenchedsf}
S_{t,\omega}:=E_{0}\big[\exp(-\int_{0}^{t}V(Z_s,\omega)ds)\big],\quad t>0,\o \in \O.
\end{equation}
Here $P_x$  is the Wiener measure on paths starting from $x\in\R^d$ and  $E_x$ is the  expectation with respect to $P_x$.

 In \cite[Theorem 0.1]{Szn94} (see also \cite[Section 5.4]{Szn98}),
 Sznitman proved a quenched large deviation principle
 for the speed of the Brownian motion in a Poissonian potential constructed from obstacles with compact support.
 Building on this work,  Armstrong and Tran \cite[Corollary 2]{ArTr14}  proved a quenched
 LDP for a wide class of Hamiltonians with stationary potentials.  
 However,  the homogenization techniques used in  \cite{ArTr14} require  some regularity of  the potential.
 In particular, the sufficient condition given for the LDP involves a finite moment of the Lipschitz norm of the potential.  The goal of this paper is to extend the quenched LDP  for the speed of the Brownian motion 
 to  stationary random potentials without regularity conditions.
 
 The sufficient conditions for this LDP  involve an integrability condition expressed in terms of the Lorentz spaces
 and the principal eigenvalue of  $-\frac 1 2 \Delta + V$. We recall these two notions before stating the LDP. 
  
The Lorentz spaces  (see for instance \cite[p.634]{Bjo10} or \cite{Boi90}) which appear in our context are defined as
$$L_\PP(d,1) =\{f: (\Omega,\mathcal{F}) \rightarrow (\mathbb{R},\mathcal{B}(\R)) \mbox{ is measurable and } ||f||_{d,1}<\infty\}$$
where $\Vert f\Vert_{d,1}=\int_0^1f^*(s)s^{(1/d)-1}ds $ and $f^*: [0,1]\rightarrow \mathbb{R}^{+}$ is the non-increasing right continuous function which has the same distribution as $|f|$. Note that $L_\PP(d,1)$ is a Banach space and
there are positive constants $c_1$ and $c_2$ such that for all $\eps>0$
\begin{equation}\label{Banach}
c_{1}\Vert f\Vert_{d}\leq\Vert f\Vert_{d,1}\leq c_{2}\Vert f\Vert_{d+\eps} \qquad \mbox{ where } \Vert f\Vert_p^{p}=\int_\Omega |f|^p d\mathbb{P}.
\end{equation}
In particular, $L^{d+\eps}_\PP \subset L_\PP(d,1)\subset L^{d}_\PP$ for all $\eps>0$.

 The principal Dirichlet eigenvalue of $-\frac 1 2 \Delta + V$   is defined as
\begin{equation}\label{Direigenvalue}
\lambda_V := \inf \left\{\int_{\R^d}\left( \frac 1 2 \vert \nabla\phii\vert^2 + V \phii^2\right) dx ; \phii \in\C_c^\infty(\R^d) , \int_{\R^d} \phii^2 dx = 1\right\}.
\end{equation}
By ergodicity, $\l_V$ is non-random.
It is closely related to the
 asymptotic behavior of the survival function.
Indeed,  \begin{equation}\label{quenchedsurvto0}
 \lim_{t\to\infty} - \frac 1 t \log S_{t,\o} =\lambda_V\qquad \PP \hbox{ - a.s.}
\end{equation} 
A proof is given  in  \cite[section 3.1]{Szn98} for non-negative potentials in the Kato class $\K_{d}^{\hbox{\small loc}}$. These include the stationary potentials which verify conditions (\ref{cond2}) and (\ref{cond22}) below.

Denote  the Lebesgue measure on $\R^d$ by $\Leb$ and
the expectation with respect to $\PP$  by $\EE$.
The Euclidean ball $\{x\in\R^d ; \vert x-y\vert < R\}$ will be denoted by $B(y, R)$
and $B(y):= B(y,1)$.

 \begin{theorem} \label{LDP}
Let $V$ be a non-negative,  stationary and  ergodic potential which verifies 
\begin{equation} \label{cond2}
\sup_{x \in B(0)}V(x,\cdot) \in L_\PP(d,1).
\end{equation}
and
  \begin{equation} \label{cond222}
   \l_V = \inf_\O V_0.
 \end{equation}

For $d=1$ or 2,  suppose moreover that there exist positive constants $ \rho, \eps$ 
and a measurable function $u : \O\to \R^d$ such that  $\PP$ - a.s.
\begin{equation} \label{cond22}
\Leb\left(\{x\in\R^d ; V(x,\o) > \eps\} \cap B(u(\o),\rho)\right) > \eps
  \mbox{ and } \mathbb{E}( | u(\cdot) |^d )<\infty.
\end{equation}
 
Then  there is a deterministic, continuous convex rate function $I : \R^d\to [0,\infty [$ 
given in (\ref{rate_function}),
with  level sets $\{x\in\R^d ; I(x)\leq c\} $ that are compact for all $c\in \R$
and such that,
$\PP$ - a.s., 
\\for all closed subsets $A$ of $\mathbb{R}^d$,
\begin{align}\label{upperLDP}
\limsup_{t\to \infty} \frac{1}{t} \log Q_{t,\omega}(Z_t \in tA) 
  \leq  -\inf_{x\in A}I(x)
\end{align}

then for all open subsets $\cO$ of $\mathbb{R}^d$,
\begin{equation}\label{lowerLDP}
\liminf_{t\to \infty}\frac{1}{t}\log Q_{t,\omega}(Z_t \in t\cO) \geq -\inf_{x\in \cO}I(x).
\end{equation}
\end{theorem}
The expression of the rate function   in terms of Lyapunov exponents allows  to prove that
 the change in regime of the Brownian motion with constant drift observed by
 Sznitman \cite[Theorem 0.3]{Szn94} 
 in a Poissonian potential associated to  obstacles with compact support actually occurs for a large class of measurable potentials.
 This phase transition was further studied by Flury \cite{Flu07, Flu08b} both in the discrete and the continuous settings.
 Concurrently, also under some regularity conditions on the potential, 
Ruess \cite{Rue14} proved the existence of the Lyapunov exponents
for Brownian motion in stationary potentials.
It  does not seem  possible to extend his results by approximating a measurable potential by
regular potentials.
Especially since Ruess  \cite{Rue16} gave an example where the Lyapunov exponents are not continuous
with respect to the potential. 

For random walks in a random potentials there is an extensive literature starting
with Varadhan \cite{Var03} who proved both a quenched and an annealed LDP
for the speed of a uniformly elliptic random walk. In his thesis, Rosenbluth \cite{Ros06}
proved a quenched large deviation principle  for a large class of random walks on $\Z^d$ with stationary transition probabilities
under an integrability condition similar to condition (\ref{cond1}).
In these works, the quenched rate function is expressed as a variational formula in terms of cocycles.

 Intensive work to extend both the class of random walks and the class of potentials for which a LDP
holds was undertaken by
  Rassoul-Agha, Sepp\"al\"ainen, Yilmaz. Some of these results, which include level-3 LDP can be found in \cite{Ras04}, \cite{Yil09}, 
 \cite{Yil11}, \cite{RaSe11},   \cite{RSY13}, \cite{RaSe14}. 
 Yilmaz and Zeitouni \cite{YiZe10} studied a class of random walks in a random environment where   the annealed
 and quenched rate functions differ. 
  
Sznitman's method, based on Lyapunov exponents, was also used to obtain
a LDP  for random walks in a random potential  in \cite{Zer98, Flu07}.
Mourrat \cite{Mou12} considered the simple random walk in an i.i.d. potential taking
values in $\lbrack 0 ; +\infty\rbrack$ and showed a LDP without assuming a moment condition on $V$.
 See also \cite{Le15}.

As a guideline for the rest of the paper, we follow \cite{Szn94}.
Along the way, we provide sufficient conditions for the intermediate results.
They are stated in section 2 and the proofs are given in section 3.

The existence of the Lyapunov exponents 
is shown in Theorem   \ref{thmlyapexp}.
for stationary potentials under a weaker integrability condition than (\ref{cond2}).
Then under (\ref{cond2}), we prove in the shape theorem \ref{shape_thm} that the convergence is uniform with respect to the direction.
The appropriate tool in this context is provided by Bj\"orklund's generalization of the shape theorem \cite{Bjo10}.
 
 The main difficulty is in the proof of (\ref{lowerLDP})
  under the additional condition on the principal Dirichlet eigenvalue.
  We will  show how key arguments of \cite[Section 2]{Szn94} can be done on a linear scale. 
  This will permit the use of the maximal inequality for cocycles  \cite[Corollary 2]{BoDe91}
  and of a technique introduced in  \cite{Boi90}.
  See also \cite{BrPi13} for an application in a different context. 
   
In the last section, we will verify  the sufficient conditions for the LDP
for  long-range Poissonian potentials. These we considered in a previous version of this paper.
They are of the form 
\begin{equation}\label{classpot}
V(x,\omega) = \sum_j W(x-\omega_j),\quad  x\in\R^d,
\end{equation}
where $\omega=(\omega_j;j\in \N) $ is a Poisson cloud in $\R^d$, $d\geq 1$ and $W(x) = \vert x \vert^{-\gamma}\wedge 1$ with 
$\g>d$. 
Moreover, for these potentials, it is possible to show
that $I(x) > 0 $ for $x\in\R^d\setminus \{0\}$.

Potentials constructed in section \ref{lacoin}
from a Boolean model also verify the sufficient conditions of the LDP.
We end the last section with the presentation a model
introduced by Ruess \cite{Rue16} which does not have decorrelation properties but still verifies a large deviation principle.

{\bf Notations.}
For $y \in \mathbb{R}^d$ and $R>0$, the Euclidean norm of $y$ is denoted by $\vert y\vert$ and
 $B(y, R)$ is the Euclidean ball $\{x\in\R^d ; \vert x-y\vert < R\}$.
 $B(y)$ stands for the unit ball $B(y,1)$ and  $H(y)=\inf \{ t\geq 0: Z_t \in \oB(y)\}$ is the hitting time of $\oB(y)$, the closure of $B(y)$.
For  an open set  $D \subset \mathbb{R}^d$,
$  T_D :=\inf \{t\geq 0, Z_t \notin D \} $  and $\C_c^\infty(D)$ is the space of infinitely differentiable functions with compact support in $D$.
For $x\in\R^d$, $[x]$ is the element of $\Z^d$ closest to $x$, with some fixed rule for ties.

The Lebesgue measure on $\R^d$ is denoted by  $\Leb$ and
the volume of the unit ball of $\R^d$   by $\L_d$.
The principal Dirichlet eigenvalue of $  - \frac1 2 \Delta$ in the unit disk is denoted by $\lambda_d$.

For a random variable $X$ and for $A\in\F$, let  $E[X, A] := E(X\1_A)$.

The constants, whose value may vary from line to line, are denoted by $c$ or $C$.
Some are numbered for subsequent reference.

\section{Main Results}\label{MResults}
\cleqn

In this section, 
the existence of Lyapunov exponents of a
 Brownian motion in a stationary potential
  and the shape theorem will be proved under appropriate moment conditions.
  Then we will show how Sznitman's method leads to a large deviation principle.
 
Recently,   Ruess  \cite{Rue14} considered  Brownian motion in a stationary potential.
Inspired by Schr\"oder \cite{Sch88},
he showed  the existence of Lyapunov exponents  for a large class of potentials and
he expressed them in terms of a variational formula.
However, the existence of Lyapunov exponents by itself follows from the subadditive theorem
 under much weaker assumptions on the potential.

For  $ x, y \in \mathbb{R}^d$ and $\o\in\O$,  define
\begin{align*}
e(x,y,\omega)  & := E_x[\exp(-\int_{0}^{H(y)}V(Z_s,\omega)ds ),H(y)<\infty]\\
a(x,y, \omega)  & :=-\log e(x,y,\omega).
\end{align*}

The measurability of $e(x,y,\o)$ can be verified by standard arguments.
 It rests on the hypothesis that $(x,\o)\mapsto \tau_x\o$\ is measurable on $\R^d\times \O$.
 Moreover, under the condition $\EE[\sup_{x\in B(0)} V(x,\cdot)] <\infty$, 
  the potential $V$ locally belongs to the  Kato class $\K_d^{\text loc}$ and
  the probabilities $e(x,y,\omega) $  are strictly positive. (cf. \cite[sections 1.2 and 5.2]{Szn98}).

We introduce the \textit{Green measure} relative to the potential $V$:
\begin{equation}\label{qgreen}
G(x,A,\omega) := E_x\big[\int_{0}^{\infty}\mathbf{1}_{A}(Z_t)\exp(-\int_{0}^{t}V(Z_s,\omega)ds)dt\big]
\end{equation}
where $x\in \mathbb{R}^d$, $\omega\in \Omega$ and $A $ is a Borel subset of $\mathbb{R}^d$. $G$ can be interpreted
as the expected occupation time measure of Brownian motion killed at rate $V(\cdot, \omega)$.
We define $g(x,y,\omega)$ as  the density function relative to the Green measure
and we call it the \textit{Green function}.
The existence of $g$ is proved in \cite[(2.2.3)]{Szn98}.

We show in the next theorem that the Green function
as well as the probabilities $e(x,y,\omega)$ have  exponential decay rates which are called Lyapunov exponents.
Theorem \ref{shape_thm} shows that, under a stronger moment condition, the convergence to the Lyapunov exponents is uniform with respect to the directions.

\begin{theorem} [Existence of Lyapunov exponents]\label{thmlyapexp}
Let $V$ be a non-negative,  stationary and  ergodic potential which verifies 
\begin{equation} \label{cond1}
\EE[\sup_{x\in B(0)} V(x, \cdot)] <\infty.
\end{equation}

For $d=1, 2$, assume moreover that (\ref{cond22}) holds. 

Then there is a non-random semi-norm $\alpha(\cdot)$ on $\mathbb{R}^d$ such that $\mathbb{P}$-a.s. and in $L^1_{\mathbb{P}}$, for all $x\in \mathbb{R}^d$
\begin{equation}\label{equa9}
\lim_{r\to \infty}\frac{1}{r}a(0,r x,\omega)=\lim_{r\to \infty}\frac{1}{r}\mathbb{E}[a(0,r x,\omega)]=\inf_{  r >0}\frac{1}{r}\mathbb{E}[a(0,r x,\omega)]=\alpha(x).
\end{equation}
$\alpha$ is called the quenched Lyapunov exponent. 

 $a(0,x,\o)$ can be replaced by $-\log g(0,x,\o)$ in (\ref{equa9}).
\end{theorem}

Bj\"orklund \cite{Bjo10} extended to a very general context the shape theorem proved in
\cite{CoDu81}  for first-passage percolation with independent passage times and in
 \cite{Boi90} for stationary passage times.
 This theorem can be applied in our framework.
 
\begin{theorem}[Shape theorem] \label{shape_thm} 
Let $V$ be a non-negative,  stationary and  ergodic potential which verifies  (\ref{cond2})  and (\ref{cond22}).

Then $\PP$ - a.s., as $x\to \infty$, $x\in\R^d$,
\begin{equation} \label{equa001}
\frac{1}{|x|}|a(0,x,\omega)-\alpha(x)|\to 0\quad 
\end{equation}
 $a(0,x,\o)$ can be replaced by $-\log g(0,x,\o)$ in (\ref{equa001}).
\end{theorem}

For the proof of theorems \ref{thmlyapexp} and \ref{shape_thm}, we need to define
\begin{align}\label{defd}
d(x,y,\omega):= \max \left(-\inf_{B(x)} \log e(\cdot, y, \omega),-\inf_{B(y)} \log e(x, \cdot, \omega)\right),\ \  x,y\in\R^d,\o\in\O.
\end{align}
By using the strong Markov property of Brownian motion, it is simple to verifiy 
that $d(\cdot,\cdot,\omega)$ is a semi-norm on $\mathbb{R}^d$. 
By \cite[Lemma 5.2.1]{Szn98}, $\PP$ - a.s. $d(\cdot,\cdot,\omega)$ defines a distance on $\R^d$
which induces the usual topology. These properties still hold in a stationary potential.

Theorem \ref{shape_thm}  is first proved for $d(0,x,\o)$.
Then lemma \ref{lemme3} and lemma \ref{lemme01} allows
 to replace $d(0,x,\o)$ by  $a(0,x,\o)$ or $-\log g(0,x,\o)$ in equation (\ref{equa001}).

We first give estimates   to compare the   quantities
$a(x,y,\omega)$,  $ -\log g(x,y,\omega)$  to  $d(x,y,\omega)$.

 Define
\begin{equation*}
F_0(\omega):= \log^{+}(\int_{B(0) \times B(0)} g(x,y,\omega)dx dy)+ \sup_{\overline{B}(0)}V(\cdot, \omega), 
\quad \o\in \Omega
\end{equation*}
and let  $F(x,\o) := F_0(\tau_x\o)$.

The proof of the following lemma can be found  in  \cite[Proposition 5.2.2]{Szn98}.
The proof is very  general as it requires only basic notions of potential theory.

\begin{lemma} \label{lemme3}
Let $V$ be a non-negative,  stationary and  ergodic potential which verifies 
 condition (\ref{cond1}) and (\ref{cond22}).
 
 Then there exists a positive constant $C$ such that for $ x,y \in \mathbb{R}^d, |x-y|>4$, $\mathbb{P}$ - a.s.
$$
\max(|d(x,y,\omega)+\log g(x,y,\omega)|, |d(x,y,\omega)-a(x,y,\omega)|) 
\leq C(1+F(x,\omega)+F(y,\omega)).
$$
\end{lemma}

\begin{lemma} \label{lemme01}
Let $V$ be a non-negative,  stationary and  ergodic potential.

(i) If (\ref{cond1}),  and (\ref{cond22}) when $d=1$ or 2, hold,  then 
 for all  $x\in \mathbb{Z}^d$, $\mathbb{P}$ - a.s.,
\begin{equation} \label{DirectionalF}
\lim_{k \to \infty}\frac{F(kx,\omega)}{k}= \lim_{k \to \infty}\mathbb{E}\frac{F(kx,\omega)}{k}=0.
\end{equation}
(ii) If (\ref{cond2}), and (\ref{cond22}) when $d=1$ or 2, hold, then
 $\mathbb{P}$ - a.s.,
\begin{equation} \label{shapeF}
\lim_{x\to \infty, x\in \mathbb{Z}^d}\frac{F(x,\omega)}{|x|}=0.
\end{equation}
\end{lemma}

The rate function of large deviation principle will be given in terms of 
the Lyapunov exponents
$\a_\l(x)$ associated with the potential $\l + V$ where  $\l \geq -\uV$ where $\uV := \inf_\O V$.
The essential properties of $\a_\l$ are gathered in the next lemma.
The upper bound (\ref{boundonalpha}) should be compared with
\cite[(5.2.31)]{Szn98} and with \cite[(65)]{Zer98} for a random walk in a random potential.

Note that, as in  \cite[section 6]{Zer98}, the results will be stated in terms of $\uV$ as it highlights the role of the principal eigenvalue $\l_V$
and it facilitates the comparison with the results from stochastic homogenization. 

\begin{lemma}\label{squarootl}
Let $V$ be a non-negative,  stationary and  ergodic potential which verifies   (\ref{cond2})  and (\ref{cond22}).  

Then $(\l,x)\mapsto \a_\l(x)$ is a continuous function on $\lbrack -\uV,\infty\lbrack \times \R^d$,

for $x\in\R^d$, $\l\mapsto\a_\l(x)$ is a concave increasing function on $\lbrack -\uV,\infty\lbrack  $,

for all $x\in\R^d$  and $\l\geq -\uV$,
\begin{align}\label{boundonalpha}
 \sqrt{2(\l+\uV})\vert x\vert \leq \alpha_\l(x) \leq \vert x\vert\sqrt {2 \left(\l  + \l_d + \EE \sup_{B(0)  } V \right)}
\end{align}
and  for all $x\neq 0$,
\begin{align}\label{squaroot0}
\a'_\l(x)_- \to 0 \hbox{   as   } \l \to \infty. 
\end{align} 
\end{lemma}

 In \cite[Proposition 2.9]{Szn98}, the lower bound for the Lyapunov exponents has the form
$$\a_\l(x) \geq \max(\sqrt{2(\l+\uV)}, C) \vert x\vert,\quad x\in\R^d$$ for some positive constant $C$.
In particular this implies the non-degeneracy of $\a_0(\cdot)$.
But the proof requires specific properties of  Poissonian potentials.
This lower bound is proved for some long-range potentials in  \cite[(2.87)]{Le15}.

 The rate function of the LDP will be given by
\begin{equation}\label{rate_function}
I(x):=\sup_{\lambda \geq -\uV} (\alpha_{\lambda}(x)- \lambda ) - \uV, \hspace{1cm}x\in \mathbb{R}^d.
\end{equation}
Bounds on the rate function are easily obtained from 
the estimates on the Lyapunov exponents given in (\ref{boundonalpha}). 
When combined with the convexity properties of the Lyapunov exponents,
we obtain the following properties of the rate function. See also \cite[Lemma 5.4.1]{Szn98}.

 \begin{lemma}\label{ratefctest}
 Let $V$ be a non-negative,  stationary and  ergodic potential which verifies  (\ref{cond2})  and (\ref{cond22}).
 
 Then for all $x\in\R^d$,
 \begin{align*}
 \frac {\vert x\vert^2} 2 \leq I(x) \leq \frac {\vert x\vert^2} 2 +\l_d + \EE \sup_{B(0)  } V 
 \end{align*}
 and  $I:\R^d\to[0,\infty[$ is a non-negative convex continuous function such that the sets $\{x\in \R^d ; I(x) \leq c\}$
 are compact for all $c\in\R$.
 \end{lemma}

 Armstrong and Tran \cite{ArTr14} obtained a large deviation principle for a diffusion 
 in a stationary convex Hamiltonian with some regularity and under a weak coercivity condition.
  In the particular case of a Brownian motion in a random potential,
  the Hamiltonian is given by
 \begin{align}\label{hamiltonian}
H(p, y) := \frac 1 2 p^2 - V(y,\o),\quad   p, y\in  \R^d.
\end{align}
Although it does not appear explicitly in \cite{Szn98},
a central object in stochastic homogenization is the effective Hamiltonian $\oH$
which appears in the homogenized problem.
It is a non-random, continuous and convex function from $\R^d$ to $\R$.
It also verifies, see \cite[section 6]{ArTr14}, 
\begin{align}\label{minV}
\oH(0) = \min_{p\in \R^d}\oH(p) = -\l_V.
\end{align}

The rate function of the LDP principle is given 
in \cite[Corollary 2]{ArTr14} by,
\begin{align*}
I_{AT}(x) := \oL(x) + \oH(0),\quad x\in\R^d,
\end{align*}
where $\oL$ is Legendre-Fenchel transform of $\oH$, that is  $\oL(x) := \sup_{p\in\R^d} (p\cdot x -\oH(p))$.

Note that the  estimates on $\oH(p)$ given in \cite[Lemma 3.1]{ArTr14}  lead to estimates on the rate function.
Therefore, as in lemma \ref{ratefctest},
the rate function of a wide class of Hamiltonians
is a non-random convex and continuous function with compact level sets.

In  \cite[(3.2)]{ArTr14},
the non-random functions $\om_\mu (x)$, which are analogous to the Lyapunov exponents,
are also  expressed in terms of
$\oH$ as $\om_\mu (x) = \sup_p \{p\cdot x ; \oH(p) \leq \mu\}$
with the convention that $\sup\emptyset = -\infty$.
 
Then to see that, under the condition (\ref{cond222}),
the  rate function $I_{AT}$ coincide with  the rate function given in (\ref{rate_function}),
one can proceed as in \cite[section 1.3]{ArTr14} :  For $x\in \R^d$, by (\ref{minV}),
\begin{align*}
\oL(x)  & := \sup_{p\in\R^d} (p\cdot x -\oH(p))\\
& = \sup_{\mu\geq\oH(0)}\sup_{p\in\R^d} \{p\cdot x -\oH(p) ; \oH(p) \leq \mu\} \\
& = \sup_{\mu\geq\oH(0)}\sup_{p\in\R^d} \{p\cdot x - \mu ; \oH(p) \leq \mu\} \\
& = \sup_{\mu\geq \oH(0)} \{\om_\mu(x) - \mu \}.
\end{align*}
With a "gauge theorem" \cite[chap 4]{ChZh95}, one could also give an analogue of $\oH_*$ by
 describing $\l_V$ in terms of the existence of a solution of $\frac {\Delta} 2 u - Vu+ \lambda u =0$
in the appropriate Sobolev space for an increasing sequence of domains.

\section{Proofs}

\textbf{Proof of lemma \ref{lemme01}}

For $d\geq 3$, 
there is a positive constant $C$ such that $\mathbb{P}$ a.s. for all  $x\neq y$,
$g(x,y,\omega) \leq   C |x-y|^{2-d} $. Hence
$$
\int_{B(0) \times B(0)} g(x,y,\omega)dx dy \leq \int_{B(0) \times B(0)} C |x -y |^{2-d}dxdy < \infty.
$$
Fix $x\in \mathbb{R}^d$. Put $X_k(\omega):= \sup_{\overline{B}(kx)}V(\cdot,\omega)$.
By condition (\ref{cond1}), $(X_k ;k\geq 0)$ is a stationary sequence of non-negative random variables with finite expectation.
Then by Borel-Cantelli lemma, $\mathbb{P}$ - a.s., $ \lim_{k\to \infty}\frac{X_k}{k}=0$.
It follows that condition (\ref{DirectionalF}) is verified for $d\geq 3$.

For $d=1$ or $2$, assume that condition (\ref{cond22}) is verified for some positive numbers  $ \rho, \eps$
and for $u:\O\to\R^d$ such that $\EE\vert u\vert < \infty$.

Consider $D :=B(0, \vert u \vert +2\rho +1)$.
Construct two increasing sequences of stopping times with respect to the natural right continuous filtration $(\mathcal{F}_t)$ on $\C(\mathbb{R}^{+}, \mathbb{R}^d)$. These stopping times describe the successive times of return to $\overline{B}(0)$ and exit times from $D$ of the Brownian motion
$$R_1:=\inf \{t \geq 0: Z_t \in \overline{B}(0) \},\qquad T_1:= \inf \{t\geq R_1, Z_t \notin D\}$$
and by induction for $n\geq 1$, 
$  R_{n+1}=R_1\circ \theta_{T_n}+T_n,\ \  T_{n+1}=T_1 \circ \theta_{R_n}+R_n$
where $\theta_t$, $t\geq 0$ is the canonical shift on $\C(\mathbb{R}^{+}, \mathbb{R}^d)$.

Since the Brownian motion is recurrent when $d=1$ or $2$, the stopping times are a.s. finite and
$$ 0\leq R_1 <T_1<R_2<T_2< \cdots <R_n <T_n\cdots \mbox{ and $R_n,T_n \uparrow \infty$ }.$$

We now have for $x\in\R^d$, 
\begin{align}\label{equa0010}
\int_{\overline{B}(0)}g(x,y,\omega)dy 
&=\int_{0}^{\infty}E_{x}[\mathbf{1}_{\overline{B}(0)}(Z_t)\exp(-\int_{0}^{t}V(Z_s)ds)]dt \notag\\
&=E_{x}[\int_{0}^{\infty}\mathbf{1}_{\overline{B}(0)}(Z_t)\exp(-\int_{0}^{t}V(Z_s)ds)dt] \notag\\
&= E_{x}[\sum_{i\geq1}\int_{R_i}^{T_i}\mathbf{1}_{\overline{B}(0)}(Z_t) \exp (-\int_{0}^{t}V(Z_s)ds)dt] \notag\\
&\leq \sum_{i=1}^{\infty}E_{x}[\exp(-\int_{0}^{R_i}V(Z_s)ds)\int_{R_i}^{T_i}\mathbf{1}_{\overline{B}(0)}(Z_t)dt] \notag\\
&= \sum_{i=1}^{\infty}E_{x}\big[\exp(-\int_{0}^{R_i}V(Z_s)ds)E_{Z_{R_i}}[\int_{0}^{T_1}\mathbf{1}_{\overline{B}(0)}(Z_t)dt]\big] \notag\\ 
& \hbox{by the strong Markov property,} \notag\\
&\leq \sup_{x\in \overline{B}(0)} E_x(T_D) \sum_{i=1}^{\infty}E_{x}[\exp(-\int_{0}^{R_i}V(Z_s)ds)]   \notag\\
&\leq C( \vert u\vert + 2\rho + 1)^2 \sum_{i=1}^{\infty}E_{x}[\exp(-\int_{0}^{R_i}V(Z_s)ds)].
\end{align}
Now, for $i\geq 1$, by the strong Markov property and by induction, for all $x\in \overline B(0)$,
\begin{align}\label{equa0011}
E_{x}[\exp(-\int_{0}^{R_{i+1}}V(Z_s)ds)]&\leq E_{x}[\exp(-\int_{0}^{R_i}V(Z_s)ds)E_{Z_{R_i}}[\exp(-\int_{0}^{T_D}V(Z_s)ds)]] \notag\\
&\leq E_{x}[\exp(-\int_{0}^{R_{i}}V(Z_s)ds)]\cdot c(\omega) \leq c(\omega)^i
\end{align}
where 
\begin{equation}\label{equa0012}
c(\omega) :=\sup_{x\in \overline{B}(0)}E_x[\exp(-\int_{0}^{T_D}V(Z_s)ds)].
\end{equation}

Note that a lower bound on the heat kernel in a region of $\R^d$ as the one 
obtained from \cite[Lemma 2.1]{Szn98}
(or more generally \cite[Theorem 3.3.5]{Dav89})
 is enough to deduce that given $\rho>0$ there is $\eta= \eta(\rho)>0$ such that for all
measurable $A\subset  B(0,\rho)$  and for all $x\in \overline B(0,\rho)$
\begin{equation}\label{vandenberg}
E_{x} [\int_{0}^{T_{B(0,2\rho)}} \1_{A}(Z_{s})ds ]> \eta \Leb(A).
\end{equation}

Now let $A:= \{V(\cdot,\o) > \eps\} \cap B(u,\rho)$ and let $Y:= \int_{0}^{T_{B(0,2\rho)}} \1_{A}(Z_{s})ds$. 
Then by (\ref{vandenberg}), (\ref{cond22}) and by Cauchy-Schwarz, there is a constant $C>0$
such that for all $x\in \overline B(u,\rho)$,
 $$\eps \eta \leq E_{x} (Y) \leq E(Y ; Y> \eps\eta/2)+ \eps\eta/2\leq (E_{x}(Y^2)P_{x}(Y>\eps\eta/2) )^{1/2}+ \eps\eta/2.$$
 Hence for all $x\in \overline B(u,\rho)$,
\begin{equation}\label{ebdavies}
P_{x} \left [\int_{0}^{T_{B(0,2\rho)} } \1_{A}(Z_{s})ds >\eps\eta/2 \right]> \left(\frac{\eps\eta} 2\right)^2 \frac 1 {E_{x}(Y^2)} 
> \frac C {\rho^4} \left(\frac{\eps\eta} 2\right)^2 . 
\end{equation}
Moreover,  by the tubular estimate \cite[p. 198]{Szn98}, there is a positive constant $C$
such that for all $t >0$ and $x\in \overline B(0)$,
\begin{align}
P_{x}\left[\sup_{0<s<t } \vert Z_s -( x_1 + \frac{s}{t} (u-x_1))\vert <  \rho\right]
  & \geq C\exp\left( -\l_d\frac{t}{\rho^2} - \frac{1}{2t}\vert u -x_1 \vert^2\right)\notag\\
  & \geq C\exp\left( -\l_d\frac{t}{\rho^2} - \frac{1}{t}(\vert u\vert^2  +1 )\right).\label{tubest}
\end{align}
Recall here that $\lambda_d$ is the principal Dirichlet eigenvalue of $  - \frac1 2 \Delta$ in the unit disk.
Hence, by taking $t=\vert u\vert+ 1 $ in (\ref{tubest}),
\begin{align}\label{tubexit2}    
P_{x}[T_D> H_{B(u,\rho)}]   
  &\geq   P_{x}\left[\sup_{0<s<t } \vert Z_s -( x + \frac{s}{t} (u-x))\vert <  \rho\right] \notag\\
  & \geq    C \exp\left( - (1+ \frac {\lambda_d} {\rho^2}) (\vert u\vert + 1 )\right).
  \end{align}

Then by (\ref{ebdavies}), (\ref{tubexit2})  and by the strong Markov property, for all   $x\in \overline B(0)$
\begin{align}\label{defpo}
P_{x}\left[ \int_{0}^{T_{D}}\textbf{1}_{A}(Z_s)ds >  \eta\eps/2 \right]
& \geq P_{x}\Big[ T_D> H_{B(u,\rho)}, P_{Z_{H_{B(u,\rho)}}} \big(\int_{0}^{T_{B(u,2\rho)}}\textbf{1}_{A}(Z_s)ds>  \eta\eps/2 \big)\Big] \notag \\
  & >    C \exp\left( - (1+ \frac {\lambda_d} {\rho^2}) (\vert u\vert + 1 )\right)  \left(\frac{\eps\eta} {2\rho^2}\right)^2 := p_0(u).
\end{align}
This provides the following upper bound for $c(\o) $ defined in (\ref{equa0012}).
\begin{align*}  
c(\o)  & = \sup_{x \in \overline B(0)} E_{x}[\exp(-\int_{0}^{T_D}V(Z_s)ds)] \\
& \leq \sup_{x \in \overline B(0)} E_{x}[\exp(-\eps \int_{0}^{T_{D}}\textbf{1}_{A}(Z_s)ds) ] \\
&\leq \exp (-\eta\eps^2 /4 )p_0(u) + 1-p_0(u) \\
 & = 1-p_0(u) (1-e^{- \eta\eps^2 /4}). 
\end{align*}
Then by  (\ref{equa0010}) and (\ref{equa0011}),  for all   $x\in \overline B(0)$,
\begin{align}\label{equa0020}
\int_{\overline{B}(0)} g(x,y,\omega) dy 
& \leq C( \vert u\vert + 2\rho + 1)^2 \sum_{i=1}^{\infty}E_{x}[ \exp(-\int_{0}^{R_i}V(Z_s)ds ) ] \notag\\
& \leq C( \vert u\vert + 2\rho + 1)^2 \frac{1}{1-c(\omega)}.  
\end{align}
Therefore by (\ref{defpo}) and (\ref{equa0020}),
\begin{align}\label{Green}
\log^+ \int_{B(0) \times B(0)} g(x,y,\omega)dx dy  & \leq C [1  +  \log^+ \vert u\vert -  \log (1 -c(\omega))]\notag\\
& \leq C [1+  \log^+ \vert u\vert -  \log p_0(u) ]\notag\\
&  \leq  C[1 +  \log^+ \vert u\vert +  \vert u\vert ].
\end{align}
Since $\EE(\log^+\vert u\vert)\leq \mathbb{E}(\vert u\vert )<\infty$, the lemma follows for $d=1,2$.

(\ref{shapeF}) follows from (\ref{Green}) and the fact that if $X(x), x\in \mathbb{Z}^d$ are identically distributed with
$\EE (\vert X(0)\vert^d )<\infty$, then $\lim_{\vert x\vert\to \infty } X(x)/|x| =0$, $\mathbb{P}$ - a.s.
  by Borel-Cantelli lemma.
\qed

\textbf{Proof of theorem \ref{thmlyapexp}}
 For a fixed $x\in \mathbb{R}^d \backslash \{0\}$, consider
$$X_{s,r}:=d(sx,rx,\omega), \hspace{1cm} 0\leq s \leq r$$
where $d$ was defined in (\ref{defd}). We have that
\begin{itemize}
\item[\rm{(i)}] $X_{s,r} \leq X_{s,u}+X_{u,r}$ for all $0\leq s \leq u \leq r$.
\item[\rm{(ii)}] $X_{s,r}\circ \tau_{u x} = X_{s+u,r+u}$ for all  $u\geq 0$.
\end{itemize}
Let $X_{[0,1]} := \sup\{X_{s,r} ; 0\leq s<r\leq 1\}$. The next step is to show that
\begin{itemize}
\item[\rm{(iii)}] $\mathbb{E}[X_{[0,1]}] < \infty$.
\end{itemize}
Let $z\in\R^d$ with $\vert z-x\vert > 1$. Then for  $\omega \in \Omega$ and $t>0$, we have that
\begin{align}\label{infezx}
e(z,x,\omega)&= E_z[\exp(-\int_{0}^{H(x)}V(Z_s,\omega) ds),H(x)<\infty] \notag \\
&\geq E_z[\exp(-\int_{0}^{t}V(Z_s,\omega) ds),\sup_{0\leq s \leq t}|Z_s-(z+\frac{s}{t}(x-z)|<1] \notag \\
&\geq P_z\big[\sup_{0\leq s \leq t}|Z_s-(z+\frac{s}{t}(x-z)|)<1\big] \exp \big(- t \sup_{y \in \C_1 (z,x)}V(y,\omega) \big)
\end{align}
where $\C_\rho (z,x) := \{y\in\R^d ; \inf_{0\leq s \leq 1 } | y - (z+ s (x-z))|< \rho\}$.

By the tubular estimate \cite[p. 198]{Szn98}, there exists a positive  constant $C$ 
such that for all $t>0$ and $\rho>0$,
\begin{equation} \label{equa4.73}
P_z\big[\sup_{0\leq s \leq t}|Z_s-(z+\frac{s}{t}(x-z))|< \rho \big]
  \geq C \exp \big(- t\frac {\lambda_d} {\rho^2}- \frac{|x-z|^2}{2t}\big).
\end{equation}

Set $t=\vert x - z\vert$.
Then by (\ref{infezx}) and (\ref{equa4.73}), there is a positive constant $C_0$ such that
\begin{align} \label{equa1000}
 & - \log e(z,x,\omega)  \leq  C_0 (\vert x- z\vert  \vee 1)+ \sup_{y \in \C_1 (z, x)} V(y,\omega)  \notag \\
\hbox{and,}\quad  & d(z,x) \leq C_0 (\vert x- z\vert  +2 )+ \sup_{y \in \C_2 (z, x)} V(y,\omega).
\end{align}
Hence, 
$\EE (X_{[0,1]} ) \leq C_0 (\vert x \vert  +2 )+ \EE [\sup_{y \in \C_2 (0, x)}V(y,\omega)]$ which is finite by (\ref{cond1}).

By  the continuous parameter subadditive theorem (see \cite[Theorem 1.5.6]{Kre85})
 and since we assumed that the
dynamical system is ergodic, there exists a constant $\alpha(x)$ such that $\PP$ - a.s. 
\begin{equation}\label{equa0001}
\lim_{r\to \infty}\frac{1}{r}d(0,r x,\omega)= \lim_{r \to \infty}\frac{1}{r}\mathbb{E}[d(0,r x,\omega)]
= \inf_{r>0}\frac{1}{r}\mathbb{E}[d(0,rx,\omega)]=\alpha(x).
\end{equation}
It is easy to check that $\alpha(\cdot)$ is a  semi-norm on $\mathbb{R}^d$. 

By  lemmas \ref{lemme3} and \ref{lemme01}, one can replace $d(0,x,\omega)$ by either one of $a(0,x,\omega)$, $-\log g(0,x,\omega)$ in (\ref{equa0001}).\qed

\textbf{Proof of Theorem \ref{shape_thm}}

By stationarity of the potential and by translation invariance of Brownian motion,
$d(x,y,\tau_z\omega)=d(x+z,y+z,\omega)$   for $z,y,z\in\R^d, \o\in\O$.
Moreover, by \cite[Lemma 5.2.1]{Szn98}, $d(\cdot,\cdot,\o)$ is a.s. a distance on $\R^d$.
Under the integrability condition (\ref{cond2}), it follows from (\ref{equa1000}) 
that $d(0,x,\o)$ is in $L_\PP(d,1)$ for all $x\in \mathbb{Z}^d$. 

Hence the conditions of the shape theorem  \cite[Theorem 1.2]{Bjo10} are verified.
Therefore, there exists a semi-norm $L$ on $\mathbb{R}^d$ such that 
\begin{equation} \label{equa4.97}
\lim_{|x| \to \infty, x\in \mathbb{Z}^d} \frac{d(0,x,\omega)-L(x)}{|x|}=0 \hspace{1cm} \mbox{a.s.}
\end{equation} 
But by Theorem \ref{thmlyapexp}, $\alpha(x)=L(x)$ for all $x\in\Z^d$ and consequently, 
for all $x\in\R^d$.

For $x\in \mathbb{R}^d$, denote by  $\hat{x}$ the nearest neighbor point in $\mathbb{Z}^d$ of $x$
(with some rule to break ties). Then, $|x-\hat{x}| < \sqrt{d}$ and for all $x\in \mathbb{R}^d\setminus B(0)$,
\begin{align} \label{equa4.99}
\frac{|d(0,x,\omega)-\alpha (x)|}{|x|} &\leq \frac{| d(0,x,\omega)-d(0,\hat{x},\omega)|}{|x|} + \frac{| d(0,\hat{x},\omega)-\alpha(\hat{x})|}{|x|}+ \frac{|\alpha(\hat{x})-\alpha(x)|}{|x|} \notag \\
&\leq \frac{|d(\hat{x},x,\omega)|}{|\hat{x}|} \cdot \frac{|\hat{x}|}{|x|}+\frac{|d(0,\hat{x},\omega)-\alpha(\hat{x})|}{|\hat{x}|}\cdot \frac{|\hat{x}|}{|x|}+ \frac{\alpha(x-\hat{x})}{|x|}
\end{align}
Consider successively the terms on the right hand side of (\ref{equa4.99})  above.
As in (\ref{equa1000}), for all $x\in\R^d$,
\begin{align*} 
d(\hat{x},x,\o) &\leq  C_0(\vert \hat x - x\vert +1)+ \sup_{y \in \C_2(\hat x,x)}V(y,\omega)\notag\\
&\leq C_0(\sqrt d  +1) + \sup_{B(\hat{x}, \sqrt{d}+3)}V(\cdot,\omega):= Y(\hat{x}).
\end{align*}
Since $(Y(\hat{x}), \hat{x} \in \mathbb{Z}^d)$ are identically distributed and  in $L^d_\PP$, by 
Borel-Cantelli lemma, 
\begin{equation*}
\lim_{x\to \infty} \frac{|d(\hat{x},x|)}{|\hat{x}|} \cdot \frac{|\hat{x}|}{|x|} \leq  \lim_{\hat{x}\to \infty, \hat{x} \in \mathbb{Z}^d} \frac{2Y(\hat{x})}{|\hat{x}|}=0\quad \PP \hbox{ - a.s.}
\end{equation*}
So the first term of (\ref{equa4.99}) converges a.s. to 0. 
From  (\ref{equa4.97}) and from (\ref{boundonalpha}) respectively,  the second and third terms converge to 0 a.s.
Hence, $\PP$ - a.s.,
\begin{equation} \label{equa16}
\lim_{|x| \to \infty, x\in \mathbb{R}^d} \frac{1}{|x|}|d(0,x,\omega)-\alpha(x)|=0
\end{equation}
By using (\ref{shapeF}) and lemma \ref{lemme3}, $d(0,x)$ can be replaced by $a(0,x)$ or $-\log g(0,x)$ in (\ref{equa16}).
 \hfill\qed
 
  \textbf{Proof of lemma \ref{squarootl}}
  
 The lower bound of (\ref{boundonalpha}) is proved as in \cite[Proposition 2.9]{Szn98}.
Let $\tV := V- \uV$. Then  for $x\in\R^d, \vert x\vert>1$,
 \begin{align*}
 e_\l(0,x,\o)  & \leq  E_0\left[\exp( -(\l + \uV) H(x)) \exp(-\int_0^{H(x)} \uV (Z_s,\o) ds), H(x)<\infty \right]\\
  & \leq \exp(-\sqrt{2(\l + \uV)}\vert x\vert)
 \end{align*}
since for a one-dimensional Brownian motion, for $\l\geq 0$ and $y\in\R$,
$E_0[ \exp (-\l H(y)) ] = \exp(-\sqrt{2\l}\vert y\vert)$.
 
 To prove the upper bound of (\ref{boundonalpha}), note that
  for $\l\geq -\uV, t>0$ and $\vert y\vert >1$,
\begin{align*}
e_\l(0,y,\o) &  =  E_0\left[\exp(-\int_0^{H(y)} (\l + V (Z_s,\o))ds), H(y)<\infty \right]\\
  &  =  E_0\left[\exp(-\int_0^{H(y)} (\l+\uV + \tV (Z_s,\o))ds), H(y)<\infty \right]\\
    & \geq  P_0[\sup_{0\leq s\leq t }\vert Z_s - \frac s t  y\vert < 1 ] \exp(-\l t -\uV t-\int_0^{t} h( \frac s t  y ,\o)ds)
    \end{align*}
  where   $h(z,\o) := \sup_{B(z)} \tV(\cdot,\o), z\in\R^d.$

Then by the tubular estimate \cite[p.198]{Szn98} and by the stationarity of $V$, 
\begin{align*}
-\EE \log   e_\l(0,y,\o) &  \leq    -\log P_0[\sup_{0\leq s\leq t }\vert Z_s - \frac s t  y\vert < 1 ]
    +( \l+\uV) t + \int_0^{t} \EE h( \frac s t  y ,\o)ds \\ 
 & \leq  C_0 + \l_dt + \frac {\vert y\vert^2}  {2t} +( \l+\uV) t + t  \EE h(0,\cdot) \\
 & =  C_0 +( \l_d +\l + \EE [\sup_{B(0)} V ] )t + \frac {\vert y\vert^2}  {2t}.
    \end{align*}
Let $y=nx$ and
$  t  =   \frac  {n\vert x\vert}{\sqrt{2(\l + \l_d  + \EE [\sup_{B(0)} V ]}} $
to obtain (\ref{boundonalpha}).

Since $\l \mapsto \a_\l(x)$ is a concave function on $[-\uV,\infty[$
 for $\l>  -\uV $,
$$\a'_\l(x)_- \leq \frac {\a_\l(x) - \a_{-\uV}(x)} { \l} \leq \frac {\a_\l(x)} { \l}.$$
 And by (\ref{boundonalpha}), $\a_\l(x) / \l\to 0$ as $\l\to\infty$.
  (\ref{squaroot0}) follows.
\qed

 \subsection{Proof of the upper estimate (\ref{upperLDP})}
 
We follow the arguments  of \cite[(4.6) of Theorem 5.4.2]{Szn98}. See also \cite[(69) of Theorem 19]{Zer98}.

First assume that $A$ is a compact subset of $\R^d$.
 For each  $t>0$, it is possible to choose $n_t$ points $x_1, x_2,\ldots, x_{n_t}$ in $A$
 such that $n_t$ grows at most polynomially in $t$ and 
$ tA \subset B_t:= \cup_{k=1}^{n_t}B(x_k).$
By definitions of $S_{t,\omega}$ and $Q_{t,\omega}$,  $\mathbb{P}$ - a.s. for all $\lambda \geq 0$,

$\exp (-\l t) S_{t,\omega} Q_{t,\omega}(Z_t \in tA)$
\begin{align}
 &= \exp (-\l t- \uV t) E_0[ \exp (-\int_{0}^{t}(V-\uV)(Z_s,\omega)ds),Z_t \in tA] \notag \\
 &\leq  \exp (-\l t- \uV t) \sum_{k=1}^{n_t} E_0[ \exp (-\int_{0}^{t}(V-\uV)(Z_s,\omega)ds),Z_t \in B(x_k)] \notag \\
 & =   \exp (- \uV t)  \sum_{k=1}^{n_t} E_0[ \exp (-\int_{0}^{t}(\lambda +V-\uV)(Z_s,\omega)ds),Z_t \in B(x_k)] \notag \\
 & \leq  \exp (- \uV t) \sum_{k=1}^{n_t} E_0[\exp (-\int_{0}^{H(x_k)}(\lambda +V-\uV)(Z_s,\omega)ds), H(x_k)<\infty] \notag \\
  &\qquad {} \text{since   }  \lambda +V-\uV \geq 0\notag\\
 & = \exp (- \uV t)  \sum_{k=1}^{n_t}e_{\lambda-\uV}(0,x_k,\omega) \leq \exp (- \uV t) n_t \max_{1\leq k \leq n_t}e_{\lambda -\uV}(0,x_k,\omega). \notag
\end{align}

Therefore for all $\lambda \geq  0$, 
by  Theorem \ref{shape_thm},     (\ref{quenchedsurvto0}) and under the assumption that $\l_V =\uV$, 
\begin{align} \label{equa10.0}
- \l -\lambda_V + \limsup_{t\to \infty} \frac{1}{t} \log Q_{t,\omega}(Z_t \in tA) &\leq -\uV -\inf_{A} \alpha_{\lambda-\uV}(x),\quad \PP - \hbox{a.s.} \notag\\
\hbox{Hence}\quad \limsup_{t\to \infty} \frac{1}{t} \log Q_{t,\omega}(Z_t \in tA) & 
  \leq  -\sup_{\lambda \geq 0}\inf_{x\in A}(\alpha_{\lambda-\uV}(x)- \lambda).
\end{align} 

To complete the proof, it remains to interchange the sup and the inf in (\ref{equa10.0}).
This is done by a classical argument (see for example \cite{DoVa75d} or \cite[p. 250]{Szn98}).
It does not require additional properties of the potential. Neither does the proof of the general case
when $A$ is a closed subset of $\R^d$ as can be seen from \cite[p. 250]{Szn98}.
\hfill\qed

\subsection{Proof of the lower estimate (\ref{lowerLDP})}

The following lemmas will be needed.

Denote by $[x,y]$, the line segment between the vertices $x,y$ of $\R^d$.
For $y\in\R^d\setminus\{0\}$, let $\Pi(y)$ be the hyperplane orthogonal to $y$ which contains $y/2$
and let $\S(y) := \{\xi\in\R^d ; \xi = z/\vert z\vert \text{   for some   } z\in \Pi(y) \text{   with   }\vert z\vert<\vert y\vert\}$.
For $\xi\in \S(y)$, denote by $[y,\xi, 0]$, the broken line from $y$ to 0 
consisting of the line segments  $[y, z]$ and   $[z, 0]$ where  $z\in \Pi(y)$ is such that $\xi= z/\vert z\vert$.
The path  integral of a measurable function $f:\R^d\to \R$ along the broken line $[y,\xi, 0]$ is denoted by
\begin{align*}
\int_{\lbrack y ,\xi ,0\rbrack} f := \int_0^{1}f(y + r  (z-y)) \vert z\vert  dr +  \int_0^{1} f(r z) \vert z\vert dr.
\end{align*}
By \cite[proposition 3]{BDD89}, if $h_0(\cdot) \in  L_\PP(d,1)$ then,
$\PP$ - a.s.,  $h(x, \o):= h_0(\tau_x\o)$ is locally in $L_{\R^d}(d,1)$, the Lorentz space over $\R^d$
with respect to the Lebesgue measure, (see  also \cite[pp. 31-32]{BoDe91}).
This in turn implies that $\PP$ - a.s., the function $\M h(y,\o)$ defined by
$$
\M h(y,\o) := \frac 1 {\sigma(\S(y))} \int_{\S(y)}d\sigma(\xi)\int_{[0,\xi, y]} h,
$$
where $\sigma(\cdot)$ denotes the Lebesgue (area) measure on the unit sphere of $\R^d$,
is continuous on $\R^d$.
Note that $\sigma(\S(y))$ does not depend on $y\in\R^d\setminus\{0\}$.

Then the arguments given in the proof of \cite[Theorem 7]{BoDe91} apply to  $ \M h(y,\o)$.
They lead to the following maximal inequality.

\begin{lemma}\label{maximal_lemma}
Let $h_0 \in L_\PP(d,1)$. Then there is a positive constant $c_3$ such that for $m>0$,
\begin{align*}
\PP\left(\sup_{\vert y\vert >1} \frac 1 {\vert y\vert}\M h(y,\o)> m \right) <  c_3 m^{-d} \Vert h_0 \Vert^d_{d,1}
\end{align*}
where $h(x, \o):= h_0(\tau_x\o)$.
\end{lemma}

\begin{lemma}\label{probacylinder}
Let $U: \R^d \to \lbrack 0, \infty \lbrack$ be  a measurable function. Then
there are positive constants $c_4$ and $c_5$ such that  
 for $y\in\R^d$, $\vert y\vert>1$, $t>0$ and $\xi\in\S(y)$,
\begin{align*}
 E_{y}\big[ \exp (-\int_{0}^{H(0)} & U(Z_s)ds), H(0) \leq t  \big] \notag\\
    \geq & c_4 \exp\left[-\l_d t- c_5 \frac{ \vert y\vert^2} {t} -  \frac t {\vert y\vert} \int_{\lbrack y ,\xi , 0\rbrack} M_2 \right]
\end{align*}
where $M_2 (z) := \sup_{x\in B(z,2)} U(x).$
\end{lemma}
{\bf Proof.}\   
It is possible to generalize the argument used in the proof of (\ref{boundonalpha}) to 
any broken line $\lbrack y ,\xi , 0\rbrack$ by combining the tubular estimates   \cite[p.198]{Szn98}
with the strong Markov property as follows. 

Let $z\in\S(y)$ be such that $z/\vert z\vert = \xi$.
Let $M_R (z) := \sup_{x\in B(z,R)} U(x)$,  $ R>0.$
 Then
\begin{align*}
E_{y}\big[ \exp (-\int_{0}^{H(0)} & U(Z_s)ds), H(0) \leq 2t  \big]   \\
 & \geq P_{y}[  \sup_{0\leq s\leq t }\vert Z_s  - y + \frac s t  (y-z) \vert < 1] \exp\left[ -\int_0^t M_1\left(y  - \frac  s t (y-z)\right)ds\right]\\
 &  \qquad{ }\ \times \inf_{z' \in B(z)}P_{z'}[  \sup_{0\leq s\leq t }\vert Z_s  - z' + \frac s t  z' \vert < 1] \exp\left[ -\int_0^t M_2\left(z  - \frac  s t z\right)ds\right]\\
 & \geq  C\exp\left[-\l_d t- \frac 1 {2t}\vert y-z\vert^2 -  \int_{0}^t M_2(y-\frac s t (y-z))ds \right]\\
 &  \qquad{ }\ \times \inf_{z' \in B(z)} C\exp\left[-\l_d t- \frac 1 {2t}\vert z' \vert^2 -  \int_{0}^t M_2(z -\frac s t z )ds \right] \\
  & \geq  C\exp\left[-2 \l_d t- \frac 1 {2t}(\vert y-z\vert^2+ \vert z' \vert^2) -  \frac t {\vert y-z\vert}\int_{[y,z]} M_2 - \frac t {\vert z\vert}\int_{[z,0]}M_2 \right]\\
  & \geq C\exp \left[ -2 \l_d t- \frac {C'} {t} \vert y\vert^2  -  \frac t {\vert z\vert}\int_{[y,\xi,0]} M_2  \right]
   \end{align*}
since $M_R \geq 0$ and by using the inequalities 
$$\frac 1 2 \leq \frac {\vert y\vert} 2 \leq\vert z\vert=\vert z-y\vert\leq \vert y\vert\hbox{   and   }\vert z' - z\vert\leq 1.$$
\hfill\qed

\begin{lemma}\label{ergothm}
Let $A $ be an event such that $\PP(A)> 1-\eps$ for some $\eps \in \rbrack 0 ,1/2\lbrack$. Let $v\in\R^d\setminus\{0\}$.

Then $\PP$ - a.s. for all    $\d>\frac{2\eps}{1-2\eps}$ and   for all sufficiently large $t$, there is $s\in\rbrack t , (1+\d)t\lbrack$ such that
 $\tau_{sv}\o \in A$.
\end{lemma} 

{\bf Proof.}\  
By the ergodic theorem, a.s. and for all $t$ sufficiently large,
\begin{align*}
1-2\eps & < \PP(A) -\eps \leq \frac 1 {(1+\d)t}\int_0^{(1+\d)t} \1_A(\tau_{sv}\o)ds
\leq \frac 1 {1+\d} + \frac 1 {(1+\d)t}\int_t^{(1+\d)t} \1_A(\tau_{sv}\o)ds
\end{align*}
Hence if $(1-2\eps)(1+\d)-1>0$ then 
  $ \int_t^{(1+\d)t} \1_A(\tau_{sv}\o)ds>0$.
\hfill\qed

 The principal Dirichlet eigenvalue of $-\frac 1 2 \Delta + V$   in the ball $B(0,R)$, $R>0$, is defined as
\begin{equation}\label{DireigenvalueR}
\l_{V,\o}(B(0,R))
:= \inf \left\{\int_{\R^d}\left( \frac 1 2 \vert \nabla\phii\vert^2 + V \phii^2\right) dx ; \phii \in\C_c^\infty(B(0,R)) , \int_{\R^d} \phii^2 dx = 1\right\}.
\end{equation}
From definitions (\ref{Direigenvalue}) and (\ref{DireigenvalueR}), it is clear that
\begin{equation}\label{uptord}
\l_{V,\o}(B(0,R)) \downarrow \l _{V}(\R^d)\ \text{   a.s.  as   }  R\to\infty.
\end{equation}
Similarly to  (\ref{quenchedsurvto0}), $\l_{V,\o}(B(0,R))$ is related to the survival time in $B(0,R)$.
We will need the following version of \cite[(3.1.17)]{Szn98} where an  $\inf_{z\in B(0)} $ appears.
The argument does not require a Harnack-type inequality.

\begin{lemma}\label{eigeninBR}
Let $V$ be a non-negative,  stationary and  ergodic potential which verifies (\ref{cond2})  and (\ref{cond22}).
Then $\PP$ - a.s. for all $R> 2$,
\begin{align*}
\liminf_t \frac 1 t \log\inf_{z\in B(0)}E_{z}\big[\exp(-\int_{0}^{t}V(Z_s,\omega)ds), T_{B(0,R)}> t\big] \geq -\l_{V,\o}(B(0,R)).
 \end{align*}
\end{lemma}
{\bf Proof.}\  
Note that a stationary ergodic potential $V$ which verifies (\ref{cond1}) and (\ref{cond22}),  also belongs to $\K_{d}^{\hbox{loc}}$
and proceed  as in
 \cite[section 3.1]{Szn98}. 
 Fix $R>0$. 
For $\eta >0$,  let $\phii \in\C_c^\infty(B(0,R))$, $\phii\geq 0$, be such that 
$$\l_{V, \o}(B(0,R))\leq \int_{\R^d}\left( \frac 1 2 \vert \nabla\phii\vert^2 + V \phii^2\right) dx \leq \l_{V,\o}(B(0,R))+\eta.$$
Let $r_R(t,x,y,\o)$ be the transition density of the Brownian motion in the potential $V(\cdot,\o)- \uV$ killed when exiting $B(0,R)$.
 Then  for $z\in B(0,R),$
\begin{align*}
\exp(-t(\l_{V,\o}(B(0,R))+\eta)) & \leq \iint \phii (x) e^{-t\uV}\ r_R(t,x,y)\phii (y) dxdy,\quad\hbox{by Jensen's inequality,}\\
 & \leq e^{\uV} \iint \phii (x)e^{-(t+1)\uV}\   \frac {r_R(1,z,x)} {\inf_{\supp \phii } r_R(1,z,\cdot)}  r_R(t,x,y)\phii (y) dxdy
\end{align*}
\begin{align*}
& \leq \frac{\Vert \phii \Vert_\infty^2}{\inf_{\supp \phii } r_R(1,z,\cdot)}e^{\uV} \  E_{z}\big[\exp(-\int_{0}^{t+1}V(Z_s,\omega)ds), T_{B(0,R)}> t+1\big] \\ 
 & \leq \frac{\Vert \phii \Vert_\infty^2}{\sup_{z\in B(0)}\inf_{\supp \phii } r_R(1, z ,\cdot)}e^{\uV} \   \inf_{z\in B(0)} E_{z}\big[\exp(-\int_{0}^{t+1}V(Z_s,\omega)ds), T_{B(0,R)}> t+1\big].
\end{align*}
And the result follows.
 \hfill\qed

A close examination of the proof of the following key lemma from \cite{Szn94} or  \cite{Szn98},
shows that it holds for stationary potentials under the moment condition (\ref{cond1}).

For    $v\in \mathbb{R}^d \setminus \{0\}$, $0<s_1< s_2<\infty$, $0\leq m \leq n$, define
\begin{align*}
  S_{m,n,v, s_1}  & :=H(nv)\circ \theta_{(n-m) s_1}+ (n-m) s_1,  \\
 A_{m,n,v,s_1,s_2} &  :=\{S_{m,n,v,s_1}<(n-m)s_2\}   
\end{align*}
where $\theta_t$, $t\geq 0$ is the canonical shift on $\C(\mathbb{R}^{+}, \mathbb{R}^d)$.
Note that $S_{m,n,v,s_1}$ is a stopping time and $A_{m,n,v,s_1,s_2}$ is the event that $Z_{\cdot}$ enters $B(nv)$ in the time interval $[(n-m)s_1,(n-m)s_2]$.
Consider
$$
b_{\lambda}(m,n,v,s_1,s_2,\omega):=-\inf_{z\in B(mv)}\log E_{z}\big[ \exp \big(-\int_{0}^{S_{m,n,v,s_1}}(\lambda+V)(Z_s,\omega)ds \big),A_{m,n,v,s_1,s_2}\big].
$$
 The strong Markov property implies 
 that $\{b_{\lambda}(m,n,v,s_1,s_2,\omega)\}_{m\geq 0,n\geq 0}$ is a subadditive sequence.
 A calculation similar to (\ref{equa1000})
shows that $\mathbb{E}b_{\lambda}(0,1,v,s_1,s_2)<\infty$.

\begin{lemma}[\cite{Szn94}, Lemma 5.4.3]\label{lemma_first_term}
Let $V$ be a non-negative,  stationary and  ergodic potential  which verifies   (\ref{cond2})  and (\ref{cond22}).
Then for $v\in \mathbb{R}^d, v\neq 0$, $\lambda>0$, $0<s_1 <s_2<\infty$,
\begin{equation}\label{defkappa}
\lim_{n\to\infty}\frac{b_{\lambda}(0,n,v,s_1,s_2,\omega)}{n}=\lim_{n\to\infty}\mathbb{E}\frac{b_{\lambda}(0,n,v,s_1,s_2,\omega)}{n} :=\kappa_{\lambda}(v,s_1,s_2) \in [0,\infty)
\end{equation}
Moreover, if $\lambda>0$ and 
$\rbrack s_1,s_2\lbrack \cap [\alpha'_{\lambda}(v)_{+},\alpha'_{\lambda}(v)_{-}] \neq \emptyset $ then 
\begin{equation}\label{kappaalpha}
\kappa_{\lambda}(v, s_1, s_2) \leq \alpha_{\lambda}(v).
\end{equation}
\end{lemma}
Here $\alpha'_{\lambda}(v)_{+},\alpha'_{\lambda}(v)_{-}$ are respectively the right and left derivatives of $\alpha_{\lambda}(v)$. 
To prove  (\ref{kappaalpha}), use (\ref{squaroot0}) and proceed as in  \cite[Lemma 5.4.3]{Szn98}.
(\ref{defkappa}) follows from Kingman's subadditive ergodic theorem. 

Note that for all $0\leq s_1 < s_2 <\infty$ and $\l>0$, $\kappa_{\lambda}(v, s_1, s_2) \geq \alpha_{\lambda}(v)$.

 \textbf{Proof of   (\ref{lowerLDP})}
 
 Let $\tV := V- \uV$ and denote the corresponding Lyapunov exponents by $\ta_\l(x)$. Then  $\ta_\l(x)= \a_{\l-\uV}(x)$ and under the assumption that $\l_V = \uV$, we have that 
 \begin{align}\label{lambdaeqzero}
 \l_{\tV} = 0 .
 \end{align}
 Since for any open set $\cO\subset \R^d$,
  \begin{align*}
  S_{t,\omega} Q_{t,\omega}(Z_t \in  t\cO)
 = \exp (- \uV t) E_0[\exp (-\int_{0}^{t}(V-\uV)(Z_s,\omega)ds), Z_t \in  t\cO ]
\end{align*} and 
 by   the continuity of $I(\cdot)$,
to obtain (\ref{lowerLDP}),
 it is sufficient to show that  for all  $v\in\Q^d\setminus \{0\}$ and $r>0$,
\begin{equation*} 
\liminf_{t\to \infty} \frac{1}{t} \log E_{0}\big[ \exp (-\int_{0}^{t}\tV(Z_s,\omega)ds), Z_t \in t B(v,r)\big] \geq -I(v) \quad  \mathbb{P}\mbox{ - a.s.}
\end{equation*}
where $I(v)= \sup_{\l\geq 0}(\ta_\l(v) - \l)$ as in (\ref{rate_function}).

We will need the following events.
For  positive numbers $t_0, R, \eps$, let $\A_1(t_0,R,\eps)$  be the event
\begin{align*}
 \left\{ \text{for all   } t>t_0, \frac 1 t \log \inf_{z\in B(0)} E_{z}  \big[ \exp (-\int_{0}^{t}\tV(Z_s,\omega)ds),  t< T_{B(0,R)} \big] > -\l_\tV-\eps\right\}
 \end{align*}
 \begin{align*}
 \hbox{and  for $m>0$, let   }   \A_2(m)  := \bigg\{ \text{for all   $y$, $\vert y\vert >1$,  there is a broken line $[0,\xi, y]$ from 0 to $y$} \\
   \text{with $\xi\in\S(y)$  and  such that $\int_{[0,\xi,  y]} M_2 < m \vert y\vert$} \bigg\}
\end{align*}
where  $ M_2(x,\o) := \sup_{B(x,2)} \tV(\cdot,\o)  $ for $x\in\R^d$.

Finally, let $\O'$ be the event of probability 1 where (\ref{defkappa}) holds.

For the moment, assume that for some positive numbers
$\eps' \in \rbrack0,1/4\lbrack, t_0, R, \eps$ and $m$,
\begin{align}\label{a1a2o}
\PP\left(\A_1(t_0, R,\eps)\cap \A_2(m) \cap \O' \right)> 1-\eps'.
\end{align}

Let $\o\in \A_1(t_0, R,\eps)\cap \A_2(m)\cap \O'  $. If  $4 \eps'< \d < r $ then $\d> 2\eps'/(1-2\eps')$ and
by lemma \ref{ergothm},  for all    $t$ sufficiently large
 there is $y_t \in\R^d$ such that
 \begin{align}\label{prelimi}
 3 < \vert [t]v - y_t\vert< \d t \vert v\vert \text{   and   } \tau_{y_t}\o \in \A_1(t_0, R,\eps)\cap \A_2(m).
 \end{align}

 Moreover, let $0< s _1< s_2 < 1$ and  $\l>0$ be such that 
$\rbrack s_1,s_2\lbrack \cap [\ta'_{\lambda}(v)_{+},\ta'_{\l}(v)_{-}] \neq \emptyset$.

Then for $\eta \in \rbrack 0 , 1- s_1\lbrack$ and for all  $t$ sufficiently large so that $R+\d t< rt$,
\begin{align}\label{estimation_lower_bound3}
E_{0} & \big[ \exp (-\int_{0}^{t}\tV(Z_s,\omega)ds), Z_t \in t B(v,r)\big] \geq \notag \\
 &  E_{0}\big[\exp (-\int_{0}^{S_{0,[t] ,v,s_1}}\tV(Z_s,\omega)ds) , A_{0,[t],v,s_1, s_2}\big]  \notag\\
 & \times \inf_{y\in B([t] v)}E_{y}\big[ \exp (-\int_{0}^{H(y_t)}\tV(Z_s,\omega)ds), H(y_t) \leq \eta t  \big] \notag\\
  & \times \inf_{z\in B(y_t)}E_{z}\big[ \exp (-\int_{0}^{(1-s_1)t} \tV(Z_s,\omega)ds), T_{B(y_t,R)} > (1-s_1)t \big].
   \end{align}
 
By lemma \ref{lemma_first_term}, since $\o\in\O'$ and  $  S_{0,[t],v,s_1}\geq [t] s_1$,
\begin{align}\label{first}
\liminf \frac 1 t \log  E_{0}\big[\exp (-\int_{0}^{S_{0,[t],v,s_1}}\tV(Z_s,\omega)ds),  &  A_{0,[t],v,s_1, s_2} \big] \notag \\
   \geq  & -\widetilde\kappa_\l(v) + \l  s_1  \geq -\ta_\l(v) + \l s_1.
 \end{align}
 
 By lemma \ref{probacylinder} and since $\tau_{y_t}\o \in\A_2(m)$, for some $\xi\in \S([t]v-y_t)$,
\begin{align*}
     \inf_{y\in B([t]v)} E_{y} &\big[ \exp (-\int_{0}^{H(y_t)}\tV(Z_s,\omega)ds), H(y_t) \leq \eta t  \big] \\
 & \geq    c_4\exp[ -\l_d \eta t  - \frac {c_5}  {\eta t} \vert [t]v - y_t\vert^2 
    -   \frac {\eta t } { \vert [t] v - y_t\vert}\int_{y_t + [[t]v - y_t ,\xi, 0]} M_2]  \\
     & \geq  c_4\exp[ -\l_d  \eta t- \frac {c_5} {\eta t} \vert [t]v -  y_t\vert^2 -   \eta t m  ] 
 \end{align*}
Hence
\begin{align}\label{second}
\liminf \frac 1 t  & \log \inf_{y\in B([t]v)} E_{y}\big[ \exp (-\int_{0}^{H(y_t)}\tV(Z_s,\omega)ds), H(y_t) \leq \eta t  \big]\notag \\
  & \geq 0 - \l_d \eta - c_6 \frac {\d^2} \eta \vert v\vert^2 - \eta m.
\end{align}

For the third term, since $\tau_{y_t}\o\in \A_1(t_0, R,\eps)$, by (\ref{lambdaeqzero}), whenever $(1-s_1)t>t_0$,
\begin{align}\label{third}
   \log \inf_{z\in B(y_t)}  E_{z}\big[ \exp (-\int_{0}^{T_{B(y_t,R)}}&  \tV(Z_s,\omega)ds), T_{B(y_t,R)} > (1-s_1)t \big] \notag\\
  & \geq -(1-s_1) (\l_\tV + \eps)t =  -(1-s_1) t \eps.
\end{align}
Putting together equations (\ref{estimation_lower_bound3}) - (\ref{third}), we find that
on $\A_1(t_0, R,\eps)\cap \A_2(m) \cap \O' $,
\begin{align}\label{lambdas1}
  \liminf_t \frac 1 t \log E_{0}  \big[ &\exp (-\int_{0}^{t}\tV(Z_s,\omega)ds), Z_t \in t B(v,r)\big] \notag \\
  & \geq  -\ta_\l(v) + \l s_1 - \l_d \eta - c_6 \frac {\d^2} \eta \vert v\vert^2 - \eta m -(1-s_1)  \eps.
  \end{align}
The proof will now be completed by contradiction. 
Assume that for some $v\in\Q^d\setminus \{0\}$ and for some positive numbers
$r, \eps_0$ and $\eps_1$, on an event of probability greater than $\eps_0$
\begin{align}\label{contradic}
\liminf \frac 1 t   \log E_{0}  \big[ \exp (-\int_{0}^{t}\tV(Z_s,\omega)ds), Z_t \in t B(v,r)\big] < - I(v) - \eps_1.
\end{align}
Set
\begin{align}\label{choices}
\eps'=\eps'(m):= 2 c_3m^{-d}\Vert h_0\Vert_{d,1}^d,\  \eta=\eta(m):= \frac{\eps_1}{20 m}\ \hbox{  and  }\ 
\delta=\delta(m) := 5\eps'.
\end{align}
Then for all $m$ positive, $ \eta m < \eps_1/10$ and $4\eps' < \d$.
Now, choose $m$ sufficiently large so that
\begin{align} \label{choicem}
\eps' < \min \{1/4, \eps_0\},\    \eta < \min  \{\eps_1/10\l_d, 1-s_1\},\   \d< r \  \hbox{  and  }\  c_6\frac{ \d^2} \eta v^2 < \eps_1/10 .
\end{align}
Furthermore, choose $\eps < \eps_1/10$. By  (\ref{uptord}) and by lemma \ref{eigeninBR},
take $t_0, R>0$  large enough so that 
$$\PP\left(\A_1(t_0, R,\eps)  \right)> 1-\eps'/2$$
and note that by lemma \ref{maximal_lemma}, for the choice of $m$ made in (\ref{choicem}),
$\PP\left(\A_2(m)  \right)> 1-\eps'/2 .$ 

Hence (\ref{a1a2o}) holds and for $\o \in \A_1(t_0, R,\eps) \cap \A_2(m)\cap \O' $, by (\ref{lambdas1}) and (\ref{choices}),
\begin{align}\label{lambdas11}
  \liminf_t \frac 1 t \log E_{0}  \big[ &\exp (-\int_{0}^{t}\tV(Z_s,\omega)ds), Z_t \in t B(v,r)\big] \notag \\
     & > -\ta_\l(v) + \l s_1 - \eps_1/2.
\end{align}

To complete the proof, consider two cases according to the value of  $ \ta'_{\l=0}(v)_+  $.  

Case 1 : $ \ta'_{\l=0}(v)_+ < 1 $.  Then $ \ta'_{\l}(v)_+< 1$ for all $\l>0$. Hence,   $I(v)= \ta_0(v)$.
 
For $\l>0 $ sufficiently small, (\ref{lambdas11}) leads to
$$ 
\liminf \frac 1 t   \log E_{0}  \big[ \exp (-\int_{0}^{t}\tV(Z_s,\omega)ds), Z_t \in t B(v,r)\big] \geq -I (v)- 3 \eps_1 /4 
 $$
 in contradiction with (\ref{contradic}).
  
Case 2 : $  \ta'_{\l=0}(v)_+\geq 1 $.
 As in \cite{Szn98}, let $\l_\infty(v):= \inf\{\l\in\Q ; \l>0 \text{  and  } \ta'_{\l}(v)_+< 1 \}$.
Then $$\l_\infty>0  \hbox{  and  }  \ta'_{\l_\infty}(v)_+\leq  1\leq  \ta'_{\l_\infty}(v)_-$$ 
 since we assumed that $  \ta'_{\l=0}(v)_+\geq 1 $ and by concavity of $\ta_\cdot(v)$.
 
 If  $\ta'_{\l_\infty}(v)_+<1$,   there are values of  $s_1<s_2 < 1$  such that
 for $s_1<1$ sufficiently close to $1$,
$  \rbrack s_1,s_2\lbrack \cap [\ta'_{\l_\infty}(v)_{+},\ta'_{\l_\infty}(v)_{-}] \neq \emptyset $ and $\l_\infty s_1 >\l_\infty - \eps_1/5$.  Then by (\ref{lambdas11}),
$$ 
\liminf \frac 1 t   \log E_{0}  \big[ \exp (-\int_{0}^{t}\tV(Z_s,\omega)ds), Z_t \in t B(v,r)\big] \geq -I (v)- 3 \eps_1 /5
 $$
 in contradiction with (\ref{contradic}).
  
If $\ta'_{\l_\infty}(v)_+=1$. Then  $\ta'_{\l}(v)_+<1$ for all $\l>\l_\infty$ and
 $\ta'_{\l}(v)_+\uparrow 1$ as $\l\downarrow\l_\infty$.
Therefore,  there are values of $s_1<s_2 < 1$ and $\l  >\l_\infty$ such that
if $s_1$ is sufficiently close to 1 and $\l$ is sufficiently close to $\l_\infty$,
then $ \rbrack s_1,s_2\lbrack \cap [\ta'_{\l}(v)_{+},\ta'_{\l}(v)_{-}] \neq \emptyset$ and by  (\ref{lambdas11}),
$$ 
\liminf \frac 1 t   \log E_{0}  \big[ \exp (-\int_{0}^{t}\tV(Z_s,\omega)ds), Z_t \in t B(v,r)\big] \geq -I (v)- 3 \eps_1 /5
 $$
 in contradiction with (\ref{contradic}).
  \qed
 
\subsection{Application to the Brownian motion with constant drift}
 
From the LDP for the speed of Brownian motion,  Varadhan's lemma, one can obtain a LDP
for Brownian motion in a random potential with a constant drift as in  \cite[Theorem 4.7]{Szn98}
by verifying the additional condition (2.1.9) of \cite{DeSt89}. But since an upper gaussian
estimate suffices, it is also verified for a stationary potential.
 This in turn
 leads to the observation of a transition from a sub-ballistic to a ballistic regime
according to the strength of the drift. A similar phenomenon is proved for the random walk
in a random potential in  \cite[Theorem B (a)]{Flu07} and in \cite[Remark 1.11]{Mou12}.

For $h\in\R^d$, the quenched path measures of the Brownian motion in the random potential $V$
with drift $h$ is given by
 \begin{align*}
 Q_{t,\o}^h(A) := \frac 1 {S_{t,\o}^h} E_0\left[\exp(h\cdot Z_t - \int_0^t V(Z_s)ds) ,  A\right],\quad t>0.
 \end{align*}
 where $S_{t,\o}^h = E_0\left[\exp(h\cdot Z_t - \int_0^t V(Z_s)ds) \right]$.
   
 The transition from a sub-ballistic  to a ballistic regime appears clearly when described by the dual norm 
  $\a_\l^*(y) :=  \sup_{x\neq 0} \frac{x\cdot y}{\a_\lambda(x)}$ for $y\in\R^d$ and $\l \geq -\uV$.
  
 \begin{prop} Let $V$ be a non-negative,  stationary and  ergodic potential
  which verifies  conditions (\ref{cond2}), (\ref{cond222}) and (\ref{cond22}) of Theorem \ref{LDP}.
 
 Then for $h\in\R^d$,
 \begin{center}
 $\a_{-\uV}^*(h) \leq 1$\quad if and only if\quad $\lim_t \frac 1 t \log S_{t, \o}^h   = 0$ a.s.
 \end{center}
  Moreover, when $\a_{-\uV}^*(h) > 1$, $\lim_t \frac 1 t \log S_{t, \o}^h = \l_h>0$ where $\l =\l_h$ is the unique solution of
 $\a^*_{\l}(h)=1$.

 \end{prop}
 
 The proof of \cite[Theorem 5.4.7]{Szn98} (see also \cite[section 5]{Flu07})  holds with  minor modifications. 
 In particular, note that to justify that $\a^*_\l(h)   \to 0$ as $\l\to \infty$, the inequality  (\ref{boundonalpha}),
 $\a_\l(x) \geq \sqrt{2(\l+\uV}) \vert x\vert  $, suffices.

 \section{Examples.}\label{examples}

In this section, we present some examples of potentials which verify the sufficient conditions of Theorem \ref{LDP}.

\subsection{A Poissonian potential: Lacoin's potential}\label{lacoin}
\cleqn
In this section, we show that the potentials introduced by 
Lacoin in \cite{Lac12}, \cite{Lac12b} verify the conditions of Theorem \ref{LDP}.

Their interest stems from the fact that the relations verified by their scaling exponents
differ substantially from those established by W\"uthrich \cite{Wut98c, Wut98, Wut98b}
for  a potential of the form (\ref{classpot}) where $W$ has compact support.

These potentials are constructed from a Poisson Boolean model.
 Let $\Omega := \{\omega=(\omega_i,r_i)_{i\geq 0}, \omega_i \in \mathbb{R}^d, r_i \geq 1\}$ be a Poisson point process in $\mathbb{R}^d \times \lbrack0,\infty\lbrack$, $d\geq 1$, whose intensity measure is given by $ \Leb\times \nu$ where  $\nu$ is a probability measure on $\lbrack 0,\infty\lbrack$ which depends on a parameter $\d>0$ and is defined by
\begin{equation} \label{equa1}
\nu([r,\infty\lbrack)=r^{-\delta}, \quad r\geq 1.
\end{equation}
Note that each Poisson cloud $\omega\in \Omega$ is a locally finite subset of $\mathbb{R}^d\times \lbrack0,\infty\lbrack$. Index $(\omega_i,r_i)$ so that $( |\omega_i|, i\geq 1)$ is an increasing sequence. See   \cite[section 1.4]{MeRo96} for an alternative description of this model and  \cite{Gou08} for results on the  percolative properties of the balls $B(\o_i, r_i)$.

Given $\gamma >0$, Lacoin's potential $ V:  \mathbb{R}^d \times \Omega \longrightarrow \lbrack 0, \infty\lbrack $ is defined by
\begin{equation*}  
V(x,\omega):= \sum_{i=1}^{\infty}r_i^{-\gamma}  \mathbf{1}_{B(\o_i, r_i)}(x).
\end{equation*}
The behavior of this model  depends on the positive parameters $\d$ and $ \g$.
For $\g+\d>d$, the potential is finite a.s. and the survival functions are strictly positive.
Additional basic properties of this potential are gathered in the following lemma taken from  \cite{Le15}.
\begin{prop}\label{lemme001}
$\PP$ - a.s. $V(x,\omega)$ is finite  for every $x\in\R^d$ if and only if $\g + \d -d >0$.
 In this case, 
\begin{equation*} 
\mathbb{E}[V(0)]=\frac{\mathcal{L}_d \delta}{\gamma+\delta-d}, \quad \mathbb{V}ar[V(0)]=\frac{\mathcal{L}_d \delta}{2\gamma+\delta-d}
\end{equation*}
and for all $R>0$, $s \in \R$,
\begin{equation} \label{Lexpmom}
\mathbb{E}\big[\exp\big(s \sum_{i=1}^{\infty} r_i^{-\gamma} \mathbf{1}_{\{|\omega_i|\leq r_i+R\}}\big)\big] = \exp\big( \int_{1}^{\infty}\delta \mathcal{L}_d (r+R)^{d}r^{-\delta-1}(e^{s r^{-\gamma}}-1)dr\big)
\end{equation}
is finite and there are positive constants 
$c_{7}(d,\delta,\gamma)$, $c_{8}(d,\delta,\gamma)$ such that  for all $x\in \mathbb{R}^d$,
$\vert x\vert >1$,
\begin{equation}\label{Lcov}
c_{7} |x|^{d-\delta-2\gamma} \leq \Cov(V(0), V(x)) \leq c_{8}  |x|^{d-\delta-2\gamma}. 
\end{equation}
Moreover, the potential is ergodic.
\end{prop}

In order to show that Lacoin's potential verifies  the conditions of Theorem \ref{LDP},
we first prove a  weak independence property  similar to \cite[Lemma 6]{Fuk11}.
The method previously used in \cite[Lemma 2.6]{Le15} lead to a weaker result.

\begin{lemma}\label{indepLac} Assume that   $ \g+\d - d >0 $. 

Then there is a constant $C=C(\gamma,\delta,d)$ such that for all $\eps>0 $, for all   $R_0>1$ and 
$R> C\eps^{-1/(\gamma+\delta-d)} \vee 2R_0$,
\begin{equation}\label{indeppropL}
\PP\left(\sup_{y\in B(0,R_0)} \sum_{\o_j\notin B(0,R)} r_j^{-\gamma}  \mathbf{1}_{B(\o_j, r_j)}(y) > \eps \right)
\leq  \exp(- 2 \eps  R^{\g}).
\end{equation}
\end{lemma}

{\bf Proof.}\  
Let $R>2R_0 > R_0 >1$. Then for all $y\in B(0,R_0)$ and 
$\o_j\notin B(0,R)$,
$$\vert \o_j \vert < \vert \o_j - y \vert + \vert y\vert < \vert \o_j - y \vert  + R_0< 2\vert \o_j - y \vert .$$
Hence 
\begin{equation}\label{LineqrR}
\sup_{y\in B(0,R_0)} \sum_{\o_j\notin B(0,R)} r_j^{-\gamma} \mathbf{1}_{B(\o_j, r_j)}(y) \leq \sum_{\o_j\notin B(0,R)}r_j^{-\gamma}\mathbf{1}_{\{|\omega_j|< 2r_j \}}.
\end{equation}
Moreover, by Campbell's theorem,  for all $s>0$,
\begin{align}\label{dominLacoin}
\log\EE[\exp  (s\sum_{\o_j\notin B(0,R)}r_j^{-\gamma}\mathbf{1}_{\{|\omega_j|< 2r_j \}} ) ] 
 & = \int_{1}^{\infty}\int_{\mathbb{R}^d}[\exp(s r^{-\gamma}\mathbf{1}_{\{R< |z|< 2r\}})-1]dz \delta r^{-\delta-1}dr \notag\\
 & =  \delta \mathcal{L}_d \int_{R/2}^{\infty}((2r)^d-R^d)(\exp(s r^{-\gamma})-1)r^{-\delta -1}dr \notag\\
 & \leq  C s \exp (s 2^{\gamma} R^{-\gamma}) \int_{R/2}^{\infty} r^{d -\gamma-\d-1}dr \notag\\
 & \leq   C s \exp (s 2^{\gamma} R^{-\gamma}) R^{d-\delta-\gamma}.
\end{align}
Then by  Markov's inequality, (\ref{LineqrR}) and  (\ref{dominLacoin}) 
with $s= 4 R^{\g}$, there exists $C=C(d,\delta, \gamma)$ such that for all $\eps >0$,
for all $R_0>1$, $R>(\frac{C}{2\eps})^{\frac{1}{\gamma+\delta-d}} \vee R>2R_0$
\begin{equation*}
\mathbb{P} [\sup_{y\in B(0,R_0)} \sum_{\o_j\notin B(0,R)} r_j^{-\gamma}\mathbf{1}_{B(\o_j, r_j)}(y) > \eps]  
\leq \exp(C R^{d-\delta}-4\eps R^{\gamma})< \exp (-2\eps R^{\gamma}).
\end{equation*}
\hfill\qed

We  are now ready to verify the hypothesis of Theorem \ref{LDP}.

\begin{prop}\label{qsurvival_lacoin} Assume that  $ \g+\d - d >0 $.
Conditions (\ref{cond2}), (\ref{cond22})  and (\ref{cond222})  are verified, the Lyapunov exponents $\a(\cdot)$
is a norm and
\begin{align}\label{quenchedLac}
  \l_V =\uV = 0.
\end{align}
\end{prop}

{\bf Proof.}\  
To check that all moments of $\sup_{x\in B(0)}V(x, \cdot)$ are finite, note that
$$\sup_{B(0)}V(\cdot,\omega) \leq  \sum_{i=1}^{\infty}r_i^{-\gamma}\mathbf{1}_{\{|\omega_i|\leq r_i+1\}}$$ and by (\ref{Lexpmom}),
\begin{equation*}
\mathbb{E}[\exp(\sum_{i=1}^{\infty}r_i^{-\gamma}\mathbf{1}_{\{|\omega_i|\leq r_i+1\}})]
= \exp\big( \int_{1}^{\infty}\delta \mathcal{L}_d (r+1)^{d}r^{-\delta-1}(e^{r^{-\gamma}}-1)dr\big)<\infty
\end{equation*}
(since $e^{r^{-\gamma}}-1<2r^{-\gamma}$ when $r$ is large enough and $\gamma+\delta-d>0$).

Then for $\omega \in \Omega$, set
\begin{equation*} 
u = u(\omega):= \omega_i \mbox{ such that } 1\leq r_i < 2 \mbox{ and if there exist  } j<i: |\omega_j|<|\omega_i| \mbox{ then } |r_i| \geq 2. 
\end{equation*}
In other words, $u(\omega)=\omega_i$ where $(\omega_i,r_i)$ is the point of the Poisson cloud in the set $\mathbb{R}^d \times [1,2]$ with $|\omega_i|$ minimum.
If we choose $\epsilon=  2^{-\gamma}$ and $\rho=1$,
then for all $x\in B(u,\rho)$, we have that
$$  V(x,\omega)=  \sum_{k=1}^{\infty}r_k^{-\gamma}   \mathbf{1}_{B(\o_k, r_k)}(x) \geq   
r_i^{-\gamma} \mathbf{1}_{B(\o_i, r_i)}(x)  >   2^{-\gamma}.$$
Note that $(\Leb\times\nu)(B(0,t)\times [1,2]) =\int_{|y|<t}\int_{1}^{2}\delta r^{-\delta-1}dydr= (1-2^{-\delta})\mathcal{L}_d t^d$.
 Then
$${\PP}(|u|>t)={\PP}(\mbox{no points of the Poisson cloud are in } B(0,t)\times [1,2])
=e^{-(1-2^{-\delta})\mathcal{L}_d t^d}.$$
Therefore $\mathbb{E}(|u|^d)<\infty$. 
Conditions (\ref{cond2}) and (\ref{cond22}) are verified.

Then by Theorem \ref{thmlyapexp},
the Lyapunov exponents $\alpha(x)$ exist.
Moreover, $\a$ is a norm. Indeed, 
$$  V(x,\omega)= \sum_{i=1}^{\infty}r_i^{-\gamma} \mathbf{1}_{B(\o_i, r_i)}(x) \geq   
 \sum_{i ; (\o_i, r_i) \in \mathbb{R}^d \times [1,2]} 2^{-\gamma} \mathbf{1}_{B(\o_i, 2)}(x) $$
which is a Poissonian potential constructed from a 
non-negative bounded measurable function with compact support. Then
by \cite[Theorem 5.2.5]{Szn98}, the associated Lyapunov exponents $\alpha^1(x)$ 
is a norm. 
And since  $\alpha(x) \geq \alpha^1(x)$,    $\alpha$ is also a norm.

To  verify (\ref{cond222}), we use lemma \ref{indepLac}.

We now prove (\ref{quenchedLac}).
For $x\in \mathbb{R}^d$ and $R>1$, write $V(x)=V_1(x)+V_2(x)$ where
$$V_1(x)=\sum_{\omega_j \in B(0,2R)}r_j^{-\gamma}\mathbf{1}_{B(\o_j, r_j)}(x) \hbox{  and  }
V_2(x)=\sum_{\omega_j \notin B(0,2R)}r_j^{-\gamma}\mathbf{1}_{B(\o_j, r_j)}(x).$$
For $\eps >0$, let $R $ be large enough so that (\ref{indeppropL}) holds.
Then by the independence property of the Poisson point process,
\begin{align*}
\PP(\sup_{x\in B(0,R)} V(x) < \eps) 
  & > \PP(\sup_{x\in B(0,R)} V_1(x) < \eps/2,\sup_{x\in B(0,R)} V_2(x) < \eps/2)\\
&=\PP(\sup_{x\in B(0,R)} V_1(x) < \eps/2) \PP(\sup_{x\in B(0,R)} V_2(x) < \eps/2)\\
& >  \exp (-CR^d)(1-\exp(-\eps R^{\gamma})).
\end{align*}
The last inequality follows from lemma \ref{indepLac} and the fact that $V_1(x)=0$ if no points
of the Poisson cloud are in $B(0,2R)\times [1,\infty[$. Hence, for all $\eps>0$ and for all $R$
large enough, 
\begin{equation}\label{posprobVsmall}
\PP(\sup_{x\in B(0,R)} V(x) < \eps) > 0.
\end{equation}

Let $(\eps_{\ell} ; \ell\in\N)$ be a sequence of positive numbers such that $\eps_{\ell}\to 0$ as $\ell\to \infty$.
Then there is a sequence $R_\ell \to \infty$ such that for all $\ell\in\N$,
$  \PP(\sup_{x\in B(0,R_\ell)} V(x) < \eps_\ell) > 0.$

By ergodicity, $\PP$ - a.s.  for each $\ell$ there is $z_{\ell}=z_{\ell}(\o)\in\R^d$ such that
$  \sup_{x\in B(z_{\ell},R_{\ell})} V(x) < \eps_\ell.$

Then $\uV=0$ and (see \cite[Section 3.1]{Szn98}),  $\PP$ - a.s. for all $\ell\in\N$,
\begin{align*}
 \l_{V}  & \leq  \l_{V}(B(z_{\ell},R_{\ell}))  \\
   & =  \inf \{ \int_{B(z_{\ell},R_{\ell})} [\frac 1 2 \vert \nabla u\vert^2  + Vu^2] dx,\    \int_{B(z_{\ell},R_{\ell})} u^2 dx = 1\}  \\
   & \leq  C R_\ell^{-2}+ \eps_{\ell}.
   \end{align*}
Let $\ell\to\infty$ to conclude.
\hfill\qed
 
\subsection{A Poissonian potential with polynomial tail}\label{classical}
\cleqn

\bi
\begin{prop}\label{qsurvival_classical}
Let $V$ be a potential of the form $$V(x,\o)=\sum_j W(x-\omega_j)$$ where
 $\o=(\o_j)$ is a Poisson point process on $\R^d$, $d\geq1$ with intensity given by Lebesgue measure
 and $W: \R^d\to [ 0,\infty\lbrack $ is a measurable function, not negligible and which verifies for $\g>d$ and
 for some positive constant $c_{9}$,
 \begin{align*}
   W(x)\leq c_{9} ( \vert x\vert^{-\g}\wedge 1) \text{   for all   } x\in\R^d.
 \end{align*}
Then conditions (\ref{cond2}), (\ref{cond222}) and (\ref{cond22}) are verified and the Lyapunov exponents $\a(\cdot)$
is a norm.
\end{prop}

For $\g>d$, the survival functions are strictly positive.
Precise estimates of the asymptotic behavior of the annealed survival
function were obtained by Donsker and Varadhan for $\g> d+2$ and by Pastur \cite{Pas77} and
Fukushima \cite{Fuk11} for $d<\g<d+2$. The case where $\g=d+2$ is considered by \^Okura \cite{Oku81}
and Chen and Kulik \cite{Che12, ChKu12} worked on the case $\g\leq d$.

The potential is ergodic since it is constructed from  a Poisson point process 
(see for instance  \cite[Proposition 2.6]{MeRo96}).

{\bf Proof.}\  
Note that $\PP$ - a.s.
$\sup_{x\in [-1,1]^d} V(x,\cdot) \leq c_{9}\sum_j \tW(\omega_j) $
where  $$\tW(x) :=  \1_{(\vert x\vert \leq 1+\sqrt{d})} + (1+\sqrt{d})^\g\vert x \vert^{-\g} \1_{(\vert x\vert> 1+\sqrt{d})}, \quad x\in\R^d.$$
Hence  for $\theta>0$,
$$\int_{\R^d} (\exp(\theta \tW(x))-1)dx
< \int_{\vert x\vert \leq 1+\sqrt{d}}(e^\theta - 1)dx + e^{\theta}\int_{\vert x\vert > 1+\sqrt{d}}  (1+\sqrt{d})^\g\vert x \vert^{-\g} dx <\infty.$$
Then by Campbell's theorem  for all $\theta\in\R$,  
 $\EE [\exp(\theta \sup_{x\in [-1,1]^d} V(x))] < \infty$.
This condition also appears as \it Assumption 2\rm\  in \cite{Fuk09}.
In particular,  $\sup_{x\in [-1,1]^d} V$ has finite moments of all order and 
condition (\ref{cond2})  holds.

Since $W$ is not negligible, there are   $\eps>0$ and $\rho>0$ such that 
  $\Leb(\{W>\eps\}\cap B(0,\rho))>\eps$.
  
Then condition (\ref{cond22}) of Theorem \ref{shape_thm} is verified with $u(\omega):=\omega_1$ where
$\omega=(\omega_i)_{i\geq 1}$ is an enumeration of the points of the   
Poisson cloud so that   $|\omega_i|\leq |\omega_{i+1}|, i\geq 1$. Indeed,
\begin{align*}
\mathbb{E}(|\omega_1|^d)&=d\int_{0}^{\infty}t^{d-1}\mathbb{P}(|\omega_1|>t)dt\\
&=d\int_{0}^{\infty}t^{d-1}\mathbb{P}(\mbox{ no points of the Poisson cloud are in } B(0,t))dt\\
&=d\int_{0}^{\infty}t^{d-1} \exp (-\mathcal{L}_d t^d)dt=d/ \mathcal{L}_d<\infty.
\end{align*}

And  $  \Leb( \{V( \cdot,\omega)>\eps\} \cap B(u,\rho) ) \geq
 \Leb( \{W(\cdot - \omega_1) >\eps\} \cap B(\o_1,\rho) ) > \eps$.

 Therefore, for $\gamma>d$,  by Theorem \ref{shape_thm},   
 there is a non-random semi-norm $\alpha(\cdot)$ on $\R^d$,  such that  
\begin{equation*}
\lim_{|x|\to \infty}\frac{1}{|x|}|a(0,x,\omega)-\alpha(x)|=0\quad \mathbb{P} \hbox{ -  a.s.}
\end{equation*}

Then  by comparison with a Poissonian potential
constructed from $W$ with compact support, the general argument given in \cite[pp. 234-236]{Szn98} shows that $\alpha$ is actually a norm

  To verify that $\lambda_V = \uV =0$, it is possible to proceed as in section \ref{lacoin} with the
appropriate version of lemma \ref{indepLac} given below.
  \hfill\qed
  
\begin{lemma}\label{Fuk11lem6} Assume that  $ \g > d $. 

Then there is a positive constant $C=C(\gamma,d)$ such that for all $\eps>0 $, for all   $r>1$ and 
$  R>C\eps^{-1/(\g - d)}\vee 2r$,
\begin{equation}\label{indepclas}
\PP\left(\sup_{y\in B(0,r)} \sum_{\o_j\notin B(0,R)}  \vert y-\o_j\vert^{-\g}> \eps \right)
\leq  \exp(-  4 \eps  R^{\g}).
\end{equation}
\end{lemma}

  {\bf Remark.} 
Note that, by Campbell's theorem,  for $\vert x\vert >2$,
\begin{align*}
\Cov(V(0), V(x) )  & = \int_{\R^d} W(u) W(x-u) du\\ 
  & \leq C \int_{\R^d} (\vert u\vert^{-\g}\wedge 1)(\vert x-u\vert^{-\g}\wedge 1)du \\
  & \leq  C\int_{B(0)}\vert x-u\vert^{-\g} du  + C\int_{B(0)^c}\vert u\vert^{-\g}\vert x-u\vert^{-\g} \1_{\{\vert u\vert<\vert x-u\vert \}}  du \\
   & < C \vert x\vert^{-\g} + C \int_{B(0)^c}\vert u\vert^{-\g}\left(\frac{\vert x\vert}{2}\right)^{-\g} \1_{\{\vert u\vert<\vert x-u\vert\}}  du \\
   &  < C \vert x\vert^{-\g}.
\end{align*}
Moreover, if there is a positive constant $c_{10}$ such that 
  $W(x) \geq c_{10} (\vert x\vert^{-\g}\wedge 1)$ for all $x\in\R^d$
  then $  \Cov(V(0), V(x) ) >  C \int_{B(0)}\vert x-u\vert^{-\g} du >  C \vert x\vert^{-\g}.$
 
\subsection{Ruess' potential}\label{ruess}

Ruess \cite{Rue14} gave an example of a two-dimensional Brownian motion in a stationary potential
constructed from a planar Poissonian tessellation. A line in the plane is parametrized by its distance to
the origin, denoted by $\rho$, and the angle $\theta\in [0, \pi[$ formed by the line and the horizontal axis.
Take $\rho \in \R$ with the convention that $\rho > 0$ if the line intersects the horizontal axis on the positive side
and $\rho\leq 0$ otherwise. Then consider a Poisson point process on $\R\times [0, \pi[$ with intensity measure
given by $\nu\Leb\times \mu$ where $\nu>0$ and $\mu$ is the uniform measure on $[0, \pi [$.

Fix $m, M \in  [0, \infty[$, $m<M $  and $R \in ]0, \infty[$.
Then the potential $V(x, \o) = m$ if $x$ is at a distance less than $R$ of one of the lines of the environment and $V(x,\o)=M$ otherwise.

For these potentials, conditions (\ref{cond2}) and (\ref{cond22}) are verified.
It is also clear from (\ref{Direigenvalue}) that $\l_V = m$. Hence, $\l_V = m =\uV$.
Therefore the shape theorem \ref{shape_thm} and the LDP given in Theorem \ref{LDP} hold for this family of potentials.
The regularity of the potential, as defined in \cite{Rue14}, is not needed for these results.
However, for regularized versions of the potentials,  \cite{Rue14}  gave a variational expression for quenched free energy.
 
\bibliographystyle{plain}
\bibliography{bibliographie}

\textsc{Universit\'e de Brest, Laboratoire de Math\'ematiques de Bretagne Atlantique,
CNRS UMR 6205, Brest, France}
\\
\textit{E-mail addresses:}  \texttt{daniel.boivin@univ-brest.fr},
   \texttt{lethithuhiensp@gmail.com}

\end{document}